\documentclass[11pt]{amsart}

\usepackage{hyperref}
\usepackage{amsmath,latexsym,amssymb,amsthm,array,amsfonts}

\usepackage{mathrsfs,dsfont,bbm}

\numberwithin{equation}{section}

\newtheorem{theorem}{Theorem}[section]
\newtheorem{lemma}{Lemma}[section]
\newtheorem{corollary}{Corollary}[section]

\newtheorem*{acknowledgement}{Acknowledgement}

\begin{document}

\title[The Bouleau-Yor identity for a bi-fBm]{The Bouleau-Yor
identity for a bi-fractional Brownian motion${}^{*}$}

\footnote[0]{${}^{*}$The Project-sponsored by NSFC (No. 11171062),
and Innovation Program of Shanghai Municipal Education Commission
(No. 12ZZ063)}

\author[L. Yan, B. Gao and J. Liu]
{Litan Yan${}^{1,\dag}$, Bo Gao${}^{1}$ and Junfeng Liu${}^{2}$}

\footnote[0]{${}^{\dag}$Corresponding Author (litanyan@hotmail.com)}

\date{}



\maketitle

\begin{center}
{\footnotesize {\it ${}^1$Department of Mathematics, Donghua
University, 2999 North Renmin Rd., Songjiang, Shanghai
201620, P.R. China}\\
{\it ${}^2$School of Mathematics and Statistics, Nanjing Audit
University, 86 West Yushan Rd., Pukou, Nanjing 211815, P.R. China}}
\end{center}

\maketitle


\begin{abstract}
Let $B$ be a bi-fractional Brownian motion with indices $H\in
(0,1),K\in (0,1]$, $2HK=1$ and let ${\mathscr L}(x,t)$ be its local
time process. We construct a Banach space ${\mathscr H}$ of
measurable functions such that the quadratic covariation $[f(B),B]$
and the integral $\int_{\mathbb R}f(x){\mathscr L}(dx,t)$ exist
provided $f\in {\mathscr H}$. Moreover, the Bouleau-Yor identity
$$
\left[f(B),B\right]_t=-2^{1-K}\int_{\mathbb R}f(x){\mathscr
L}(dx,t),\qquad t\geq 0,
$$
holds for all $f\in {\mathscr H}$.
\end{abstract}

\section{Introduction}\label{sec1}

The bi-fractional Brownian motion (bi-fBm) with indices $H\in (0,1)$
and $K\in (0,1]$ is a zero mean Gaussian process $B=\left\{B_t,t\geq
0\right\}$ such that $B_0=0$ and
\begin{equation}\label{sec1-eq1.1}
E\left[B_tB_s\right]=\frac{1}{2^K}\left[
\left(t^{2H}+s^{2H}\right)^{K}-|t-s|^{2HK}\right]
\end{equation}
for all $s,t\geq 0$. Clearly, if $K=1$, the process is a fractional
Brownian motion with Hurst parameter $H$. Bi-fBm was first
introduced by Houdr\'e--Villa~\cite{Hou}. The process $B$ is
$HK$-selfsimilar but it has no stationary increments. It has
H\"older continuous paths of order $\delta<HK$ and its paths are not
differentiable. An interesting property is that the bi-fBm has
non-trivial quadratic variation equal with a constant times $t$ in
the case $2HK=1$, which is similar to this of the standard Brownian
motion. That is
$$
\left[B,B\right]_t=\lim_{\varepsilon\downarrow
0}\frac{1}{\varepsilon}\int_0^t(B_{ s+\varepsilon} -B_s)^2ds
=2^{1-K}t,\qquad t\geq 0
$$
in $L^2(\Omega)$ (for this, see Russo-Tudor~\cite{Russo-Tudor}).
This motivates us to study the quadratic covariation and related to
stochastic calculus of bi-fBm with $2HK=1$. More works for bi-fBm
can be found in Es-sebaiy--Tudor~\cite{Es-sebaiy},
Jiang-Wang~\cite{Jiang2}, Kruk {\it et al}~\cite{Kruk},
Lei-Nualart~\cite{Lei-Nualart}, Russo-Tudor~\cite{Russo-Tudor},
Tudor-Xiao~\cite{Tudor-Xiao}, Shen-Yan~\cite{Shen-yan}, Yan {\em et
al}~\cite{Yan4} and the references therein.

Let now $2HK=1$ and let $B=\left\{B_t,0\leq t\leq T\right\}$ be the
bi-fBm on ${\mathbb R}$ with indices $H$ and $K$. In order to
motivate our subject, let us first recall some known results
concerning the quadratic variation and It\^o's formula. Let $W$ be a
standard Brownian motion and let $F$ be an absolutely continuous
function with locally square integrable derivative $f$, that is,
$$
F(x)=F(0)+\int_0^xf(y)dy
$$
with $f$ being locally square integrable. F\"ollmer {\it et
al}~\cite{Follmer} introduced the following It\^o's formula:
\begin{equation}\label{sec1-eq1.2}
F(W_t)=F(0)+\int_0^tf(W_s)dW_s+\frac12\left[f(W),W\right]_t.
\end{equation}
Moreover, the result has been extended to some semimartingales and
smooth nondegenerate martingales (see Russo--Vallois~\cite{Russo2}
and Moret--Nualart~\cite{Moret}). Thus, it is natural to ask whether
the similar It\^o formula for bi-fractional Brownian motion $B$ with
$2HK=1$, more general, for finite quadratic variation process $X$
holds or not. We will consider the question. Recall that a process
$X$ is said to be of finite quadratic variation if quadratic
variation $[X,X]$ is finite. For any continuous finite quadratic
variation process $X$ and twice-differentiable function $f$, we have
(see, for example, Russo-Vallois~\cite{Russo-Vallois2})
\begin{equation}\label{sec1-eq100}
f(X_t)=f(0)+\int_0^tf'(X_s)d^{-}X_s+\frac12\left[f'(X),X\right]_t,
\end{equation}
where the integral $\int_0^tf(X_s)d^{-}X_s$ is the forward
(pathwise) integral defined by
$$
\int_0^tf(X_s)d^{-}X_s=\lim_{\varepsilon\downarrow
0}\frac1{\varepsilon}\int_0^tf(X_s)(X_{s+\varepsilon}-X_s)ds
$$
and the quadratic covariation $[f'(X),X]$ of $f'(X)$ and $X$ is
defined as
\begin{equation}\label{sec1-eq1.4}
[f'(X),X]_t=\lim_{\varepsilon\downarrow
0}\frac{1}{\varepsilon}\int_0^t\left\{f'(X_{ s+\varepsilon})
-f'(X_s)\right\}(X_{s+\varepsilon}-X_s)ds,
\end{equation}
provided the limit exists uniformly in probability. However, the
formula~\eqref{sec1-eq100} is only effective on twice-differentiable
functions. It is impossible to list here all the contributors in
previous topics. Some surveys and complete literatures could be
found in Nualart~\cite{Nua2}, Russo-Vallois~\cite{Russo-Vallois3}
and F. Russo-Tudor~\cite{Russo-Tudor}. In this paper, our aim is to
prove It\^o's formula~\eqref{sec1-eq100} holds for $X=B$ with
$2HK=1$ whatever $f\in C^2({\mathbb R})$, and obtain the relation
between the forward (pathwise) integral and the Skorohod integral of
bi-fractional Brownian motion with $2HK=1$. Though our method is
only effective on bi-fractional Brownian motion, the merit here has
been to concentration fully on fBm in order to get a stronger
statement by fully using bi-fractional Brownian motion's regularity.
In the present paper, we consider the case $2HK=1$. Our start point
is to consider the decomposition
\begin{equation}\label{sec1-eq1.5}
\begin{split}
&\frac{1}{\varepsilon}\int_0^t\left\{f(B_{
s+\varepsilon})-f(B_s)\right\}(B_{s+\varepsilon}-B_s)ds\\
&=\frac{1}{\varepsilon}\int_0^tf(B_{
s+\varepsilon})(B_{s+\varepsilon}-B_s)ds-\frac{1}{\varepsilon}
\int_0^tf(B_s)(B_{s+\varepsilon}-B_s)ds.
\end{split}
\end{equation}
By estimating the two terms of the right hand side in the
decomposition~\eqref{sec1-eq1.5}, respectively, we can construct a
Banach space ${\mathscr H}$ of measurable functions $f$ on $\mathbb
R$ such that $\|f\|_{\mathscr H}<\infty$, where
\begin{align*}
\|f\|_{\mathscr H}^2:=\int_0^T\int_{\mathbb
R}|f(x)|^2\varphi_s(x)dxds+\int_0^T\int_{\mathbb
R}|f(x)x|^2\varphi_s(x)\frac{dxds}{s}
\end{align*}
with $\varphi_s(x)=\frac{1}{\sqrt{2\pi s}}e^{-\frac{x^2}{2s}}$. We
show that the {\it quadratic covariation} $[f(B),B]_t$ exists in
$L^2$ for all $t\in [0,T]$ if $f\in {\mathscr H}$. This allows us to
write the following It\^o's formulas (F\"ollmer-Protter-Shiryayev's
formula):
\begin{align}\label{sec1-eq1.6}
F(B_t)&= F(0)+\int_0^tf(B_s)dB_s+2^{K-2} \left[f(B),B\right]_t,
\end{align}
where the integral
$\int_0^{\cdot}f(B_s)dB_s$ is the Skorohod integral, $f\in {\mathscr H}$ is left continuous with right limit and $F$ is an absolutely continuous function with $\frac{d}{dx}F=f$. This extends the formula~\eqref{sec1-eq100} for bi-fractional Brownian motion $B$ with $2HK=1$. As an application we establish the following integral:
\begin{equation}\label{sec1-eq1.7}
\int_{\mathbb {R}}f(x){\mathscr L}(dx,t),\qquad t\in [0,T],
\end{equation}
and show that the Bouleau-Yor identity
\begin{equation}\label{sec1-eq1.8}
\left[f(B),B\right]_t=-2^{1-K}\int_{\mathbb {R}}f(x){\mathscr
L}(dx,t)
\end{equation}
holds provided $f\in {\mathscr H}$, where
$$
{\mathscr L}(x,t)=\int_0^t\delta(B_s-x)ds
$$
is the local time of bi-fractional Brownian motion $B$.

For $K=1$ and $H=\frac12$, the process $B$ is classical Brownian motion $W$ and the above results first are studied by Bouleau--Yor~\cite{Bouleau} and F\"ollmer {\it et
al}~\cite{Follmer}. Moreover, these have also been extended to semimartingales by Bardina--Rovira~\cite{Bardina},
Eisenbaum~\cite{Eisen1,Eisen2}, Elworthy {\it et
al}~\cite{Elworthy}, Feng--Zhao~\cite{Feng}, Peskir~\cite{Peskir1}, Rogers--Walsh~\cite{Rogers2}, Yan--Yang~\cite{Yan3}. For $K=1$ and $H\neq\frac12$, the process $B$ is a standard fractional Brownian motion $B^H$ with Hurst index $H$. Yan {\it et al}~\cite{Yan1,Yan2} studied the integration with respect to local time of fractional Brownian motion, and the {\it weighted quadratic covariation} $[f(B^H),B^H]^{(W)}$ of $f(B^H)$ and $B^H$. These deduce the fractional It\^{o} formula for new classes of functions. For $2HK=1$ and $K\neq 1$, this process is not fractional Brownian motion, and the question has not been studied. Recently, the long-range property has become an important aspect of stochastic models in various
scientific area including hydrology, telecommunication, turbulence,
image processing and finance. It is well-known that fractional
Brownian motion is one of the best known and most widely used
processes that exhibits the long-range property, self-similarity and
stationary increments. It is a suitable generalization of classical Brownian motion. On the other hand, many authors have proposed to use more general self-similar Gaussian process and random fields as stochastic models. Such applications have raised many interesting theoretical questions about self-similar Gaussian processes and fields in general. However, contrast to the extensive studies on fractional Brownian motion, there has been little systematic investigation on other self-similar Gaussian processes. The main reason for this is the complexity of dependence structures for
self-similar Gaussian processes which does not have stationary
increments. The bi-fractional Brownian motion has properties
analogous to those of fractional Brownian motion (self-similarity,
long-range dependence, H\"older paths, the variation and the
renormalized variation). However, in comparison with fractional
Brownian motion, the bi-fractional Brownian motion has
non-stationary increments and the increments over non-overlapping
intervals are more weakly correlated and their covariance decays
polynomially as a higher rate. The above mentioned properties make
bi-fractional Brownian motion a possible candidate for models which
involve long-dependence, self-similarity and non-stationary.
Therefore, it seems interesting to study the quadratic covariation
and extension of It\^o's formula of bi-fractional Brownian motion
with $2HK=1$.

This paper is organized as follows. In Section~\ref{sec2} we present some preliminaries for bi-fractional Brownian motion. In Section~\ref{sec3-0}, we establish some technical estimates associated with bi-fBm with $2HK=1$ and it seems interesting that these inequalities arising from the method. In Section~\ref{sec3}, we will construct the Banach space ${\mathscr H}$ such that the quadratic covariation $[f(B),B]$ exists in $L^2$ for $f\in {\mathscr H}$. In section~\ref{sec4-1} our main object is to explain and prove the generalized It\^o type formula~\eqref{sec1-eq1.6}. As an application we introduce the relationship between the forward (pathwise) integral and Skorohod integral
\begin{align*}
\int_0^tf(B_s)d^{-}B_s&=\int_0^tf(B_s)dB_s
+\frac12(2^{K-1}-1)\left[f(B),B\right]_t
\end{align*}
for all $f\in {\mathscr H}$. The result weakens the hypothesis of differentiability for $f$ (see Russo-Tudor~\cite{Russo-Tudor}). In
Section~\ref{sec4} we study the integral~\eqref{sec1-eq1.7} and show
that the Bouleau-Yor identity~\eqref{sec1-eq1.8} holds.

\section{Preliminaries for bi-fractional Brownian motion}\label{sec2}

In this section, we briefly recall the definition and properties of
stochastic integral with respect to bi-fBm. As a Gaussian process,
it is possible to construct a stochastic calculus of variations with
respect to $B$. We refer to Al\'os {\it et al}~\cite{Nua1} and
Nualart~\cite{Nua2} for a complete description of stochastic
calculus with respect to Gaussian processes. Here we recall only the
basic elements of this theory (see
Es-sebaiy--Tudor~\cite{Es-sebaiy}). Throughout this paper we assume
that $2HK=1$. As we pointed out before, bi-fractional Brownian
motion (bi-fBm in short) $B= \left\{B_t,0\leq t\leq T \right\}$, on
the probability space $(\Omega, {\mathscr F},P)$ with indices $H\in
(0,1)$ and $K\in (0,1]$ is a rather special class of self-similar
Gaussian processes such that $B_0=0$ and
\begin{equation}\label{sec2-eq2.1}
E\left[B_tB_s\right]=R(t,s):=\frac{1}{2^K}\left[
\left(t^{2H}+s^{2H}\right)^{K}-|t-s|^{2HK}\right],\quad \forall
s,t\geq 0.
\end{equation}
The process is $HK$-self similar and satisfies the following
estimates (the quasi-helix property)
\begin{equation}\label{sec2-eq2.2}
2^{-K}|t-s|^{2HK}\leq E\left[\left(B_t-B_s\right)^2\right]\leq
2^{1-K}|t-s|^{2HK}.
\end{equation}
Thus, Kolmogorov's continuity criterion implies that bi-fBm is
H\^older continuous of order $\delta$ for any $\delta<HK$.

Let $\mathcal H$ be the completion of the linear space ${\mathcal
E}$ generated by the indicator functions $1_{[0,t]}, t\in [0,T]$
with respect to the inner product
$$
\langle 1_{[0,s]},1_{[0,t]} \rangle_{\mathcal H}=R(s,t).
$$
The application $\varphi\in {\mathcal E}\to B(\varphi)$ is an
isometry from ${\mathcal E}$ to the Gaussian space generated by $B$ and it can be extended to ${\mathcal H}$. For $2HK=1$ we can characterize ${\mathcal H}$ as
$$
{\mathcal H}=\{f:[0,T]\to {\mathbb R}\;|\;\|f\|_{\mathcal H}<\infty\},
$$
where
$$
\|f\|^2_{\mathcal H}=\int_0^T\int_0^Tf(t)f(s)\phi(s,t)dsdt
$$
with $\phi(s,t)=
\left(s^{2H}+t^{2H}\right)^{K-2}s^{2H-1}t^{2H-1}$. Let us denote by $\mathcal S$ the set of smooth functionals of the form
$$
F=f(B(\varphi_1),B(\varphi_2),\ldots,B(\varphi_n)),
$$
where $f\in C^{\infty}_b({\mathbb R}^n)$ and $\varphi_i\in {\mathcal H}$. The {\em Malliavin derivative} $D^{H,K}$ of a functional $F$ as above is given by
$$
D^{H,K}F=\sum_{j=1}^n\frac{\partial f}{\partial
x_j}(B(\varphi_1),B(\varphi_2), \ldots,B(\varphi_n))\varphi_j.
$$
The derivative operator $D^{H,K}$ is then a closable operator from $L^2(\Omega)$ into $L^2(\Omega;{\mathcal H})$. We denote by ${\mathbb D}^{1,2}$ the closure of ${\mathcal S}$ with respect to the norm
$$
\|F\|_{1,2}:=\sqrt{E|F|^2+E\|D^{H,K}F\|^2_{{\mathcal H}}}.
$$

The {\it divergence integral} $\delta^{H,K}$ is the adjoint of derivative operator $D^{H,K}$. That is, we say that a random variable $u$ in $L^2(\Omega;{\mathcal H})$ belongs to the domain of the divergence operator $\delta^{H,K}$, denoted by ${\rm {Dom}}(\delta^{H,K})$, if
$$
E\left|\langle D^{H,K}F,u\rangle_{\mathcal H}\right|\leq
c\|F\|_{L^2(\Omega)}
$$
for every $F\in {\mathbb D}^{1,2}$, where $c$ is a constant depending only on $u$. In this case $\delta^{H,K}(u)$ is defined by the duality relationship
\begin{equation}
E\left[F\delta^{H,K}(u)\right]=E\langle D^{H,K}F,u\rangle_{\mathcal H}
\end{equation}
for any $F\in {\mathbb D}^{1,2}$. We have ${\mathbb D}^{1,2}\subset
{\rm {Dom}}(\delta^{H,K})$ and for any $u\in {\mathbb D}^{1,2}$
\begin{align*}
E\left[\delta^{H,K}(u)^2\right]&=E\|u\|^2_{\mathcal H}+E\langle
D^{H,K}u,(D^{H,K}u)^{*}\rangle_{{\mathcal H}\otimes {\mathcal
H}}\\
&=E\|u\|^2_{\mathcal H}+E\int_{[0,T]^4}D_\xi^{H,K}u_rD_\eta^{H,K}
u_s\phi(\eta,r)\phi(\xi,s)dsdrd\xi d\eta,
\end{align*}
where $(D^{H,K}u)^{*}$ is the adjoint of $D^{H,K}u$ in the Hilbert space ${\mathcal H}\otimes{\mathcal H}$. We will denote
$$
\delta^{H,K}(u)=\int_0^Tu_sdB_s
$$
for an adapted process $u$, and it is called Skorohod integral.

\begin{theorem}[It\^o's formula~\cite{Es-sebaiy}]\label{theorem-Ito}
Let $f\in C^{2}({\mathbb R})$ such that
\begin{equation}\label{sec2-Ito-con1}
\max\left\{|f(x)|,|f'(x)|,|f''(x)|\right\}\leq \kappa e^{\beta x^2},
\end{equation}
where $\kappa$ and $\beta$ are positive constants with
$\beta<(4T)^{-1}$. Suppose that $2HK=1$, then we have
\begin{align*}
f\left(B_t\right)=f(0)&+\int_0^t\frac{d}{dx}f(B_s)dB_s
+\frac12\int_0^t\frac{d^2}{dx^2}f(B_s)ds.
\end{align*}
\end{theorem}
Recall that bi-fBm $B$ has a local time ${\mathscr L}(x,t)$
continuous in $(x,t)\in {\mathbb R}\times [0,\infty)$ which
satisfies the occupation formula (see Geman-Horowitz~\cite{Geman})
\begin{equation}\label{sec2-1-eq1}
\int_0^t\psi(B_s,s)ds=\int_{\mathbb R}dx\int_0^t\psi(x,s) {\mathscr
L}(x,ds)
\end{equation}
for every continuous and bounded function $\psi(x,t):{\mathbb
R}\times {\mathbb R}_{+}\rightarrow {\mathbb R}$ and any $t\geqslant
0$, and such that
$$
{\mathscr L}(x,t)=\int_0^t\delta(B_s-x)ds=\lim_{\epsilon \downarrow
0}\frac{1}{2\epsilon}\lambda\big(s\in[0,t],|B_s-x|<\epsilon\big),
$$
where $\lambda$ denotes Lebesgue measure and $\delta$ is the Dirac
delta function. Moreover ${\mathscr L}$ has a compact support in $x$
for all $t\geq 0$ and the following Tanaka formula holds:
\begin{equation}\label{eq2.4}
|B_t-x|=|x|+\int_0^t{\rm sign}(B_s-x)dB_s+{\mathscr L}(x,t).
\end{equation}
For these see Es-sebaiy--Tudor~\cite{Es-sebaiy} and
Tudor--Xiao~\cite{Tudor-Xiao}.

\section{Some estimates associated with bi-fBm with $2HK=1$}
\label{sec3-0}
In this section we will establish some technical estimates
associated with bi-fBm. For simplicity throughout this paper we let
$C$ stand for a positive constant depending only on the subscripts
and its value may be different in different appearance, and this
assumption is also adaptable to $c$.
\begin{lemma}\label{lem2.2}
Let $2HK=1$, and for all $s,r\in [0,T],\;s\geq r$ we denote
$$
\rho^2_{s,r}=sr-\mu^2
$$
where $\mu_{s,r}=E(B_sB_r)$. Then we have
\begin{equation}\label{sec2-1-eq2.7}
r(s-r)\leq \rho^2_{s,r}\leq (1+2^{1-2K})r(s-r).
\end{equation}
\end{lemma}

By the local nondeterminacy of bi-fBm we can prove the lemma.
Moreover, one can also obtain the estimates by considering the
asymptotic property of some functions. Here, we shall prove these
estimates~\eqref{sec2-1-eq2.7} by an elementary method, and it seems
interesting that these inequalities arising from the method. We shall use the following inequalities:
\begin{align}\label{sec2-1-eq2.4}
(1+x)^\alpha&\leq 1+(2^\alpha-1)x^\alpha\\
\label{sec2-1-eq2.4--0} (1+x)^\beta&\geq 1+(2^\beta-1)x^\beta
\end{align}
for all $0\leq x\leq 1$, $0\leq \alpha\leq 1$ and $\beta\geq 1$. The
inequalities above are two calculus exercises, and they are stronger
than the well known inequalities
\begin{align*}
(1+x)^\alpha &\leq 1+\alpha x^\alpha\leq 1+x^\alpha\\
(1+x)^\beta &\geq 1+x^\beta
\end{align*}
because of $2^\alpha-1\leq \alpha$ and $2^\beta-1\geq 1$ for all
$0\leq \alpha\leq 1$ and $\beta\geq 1$. Furthermore, by applying the
inequality~\eqref{sec2-1-eq2.4} one can improve the left estimate
in~\eqref{sec2-eq2.2} as (see Yan {\em et al}~\cite{Yan4})
$$
|t-s|^{2HK}\leq E\left[(B_t-B_s)^2\right],
$$
for all $H\in (0,1)$ and $K\in (0,1]$.

\begin{proof}[Proof of Lemma~\ref{lem2.2}]
Clearly, by~\eqref{sec2-1-eq2.4} we have
\begin{equation}\label{sec2-1-eq2.8}
\left(s^{2H}+r^{2H}\right)^{K}\leq
s^{2HK}+(2^K-1)r^{2HK}=(s-r)+2^Kr.
\end{equation}
It follows that
\begin{align*}
\rho_{s,r}&=sr-\mu^2_{s,r} =sr-\frac1{4^K}\left[(s^{2H}+r^{2H})^K-(s-r)\right]^2\\
&\geq sr-r^2=r(s-r).
\end{align*}
In order to show that the right estimate in~\eqref{sec2-1-eq2.7}, we
have $\frac12\leq K\leq 1$ and
$$
(s^{2H}+r^{2H})^{2K}\geq
s^{4HK}+(2^{2K}-1)r^{4HK}=s^2+(2^{2K}-1)r^2,
$$
which deduces
\begin{align*}
\rho_{s,r}&=sr-\mu^2_{s,r}=sr-\frac1{4^K}\left[(s^{2H}+r^{2H})^K -(s-r)\right]^2\\
&=\frac1{4^K}\left(4^Ksr-(s^{2H}+r^{2H})^{2K}+2(s^{2H}
+r^{2H})^K(s-r)-(s-r)^2\right)\\
&\leq \frac1{4^K}\left\{4^Ksr-\left(s^2+(2^{2K}-1)r^2\right)+2(s
+r)(s-r)-(s-r)^2\right\}\\
&=\frac1{4^K}\left\{(4^K+2)sr-(4^K+2)r^2\right\}\\
&= (1+2^{1-2K})r(s-r)
\end{align*}
by~\eqref{sec2-1-eq2.4--0} and the inequality
$$
(s^{2H} +r^{2H})^K\leq s+r
$$
for all $s,r\geq 0$. This completes the proof.
\end{proof}
\begin{lemma}\label{lem2.1}
Let $2HK=1$. Then we have
\begin{equation}\label{sec2-1-eq2.5}
c_{H,K}\frac{r}{s}(s-r)\leq r-\mu_{s,r}\leq C_{H,K}\frac{r}{s}(s-r)
\end{equation}
and
\begin{equation}\label{sec2-1-eq2.6}
s-r\leq s-\mu_{s,r}\leq C_{H,K}(s-r)
\end{equation}
for all $s>r\geq 0$.
\end{lemma}
\begin{proof}
In order to show that the estimates~\eqref{sec2-1-eq2.5}, we have
\begin{align*}
r-\mu_{s,r}&=r-\frac{1}{2^{K}}\left[(s^{2H}+r^{2H})^K-(s-r)\right]\\
&=\frac1{2^K}s\left\{2^Kx-(1+x^{2H})^K+1-x\right\}\\
&=\frac1{2^K}s\left\{1+(2^K-1)x-(1+x^{2H})^K\right\}>0
\end{align*}
for all $x=\frac{r}{s}\in (0,1)$ by the inequality~\eqref{sec2-1-eq2.4}.
Elementary calculus can show that
\begin{align}\label{sec3-eq000-1}
\lim_{x\to 0}\frac{1+(2^K-1)x-(1+x^{2H})^K}{x(1-x)}&=2^K-1\\
\label{sec3-eq000-2} \lim_{x\to
1}\frac{1+(2^K-1)x-(1+x^{2H})^K}{x(1-x)}&=1-2^{K-1},
\end{align}
which deduce the estimates~\eqref{sec2-1-eq2.5} by continuity.

On the other hand, by the inequality~\eqref{sec2-1-eq2.4} we have
\begin{align*}
s-\mu_{s,r}&=s-\frac{1}{2^{K}}\left[(s^{2H}+r^{2H})^K-(s-r)\right]\\
&\geq s-2^{-K}\left(s+(2^K-1)r\right)+2^{-K}(s-r)\\
&=s-r
\end{align*}
for all $s\geq r$, and moreover we have
\begin{align*}
s-\mu_{s,r}&=s-r+(r-\mu_{s,r})\leq C_{H,K}(s-r)
\end{align*}
by the right estimates in~\eqref{sec2-1-eq2.5}.
\end{proof}

\begin{lemma}\label{lem2.3}
Let $2HK=1$. Then we have
\begin{equation}\label{sec2-1-eq2.9}
|E(B_{t} -B_{s})(B_{t'}-B_{s'})|\leq
(2H-1)2^{-K}\frac{(t-s)(t'-s')}{s}
\end{equation}
holds for all $T\geq t>s\geq t'>s'>0$.
\end{lemma}
\begin{proof}
For $y>0$ we define the function $x\mapsto G_y(x)$ on $[0,T]$ by
$$
G_y(x)=\left(y^{2H}+x^{2H}\right)^{K-1}.
$$
Thanks to mean value theorem, we see that there is an $\xi_l\in (s',t')$ such that
$$
G_l(t')-G_l(s')=2H(K-1)\xi_l^{2H-1}(t'-s')\left(l^{2H}+\xi_l^{2H} \right)^{K-2}.
$$
It follows from the duality relationship that
\begin{align*}
-E(B_{t}&-B_{s})(B_{t'}-B_{s'})
=-\int_s^t\int_{s'}^{t'}\frac{\partial^2}{\partial r\partial l}
R(r,l)drdl\\
&=2H(1-K)2^{-K}\int_s^t\int_{s'}^{t'}
\left(r^{2H}+l^{2H}\right)^{K-2}r^{2H-1}l^{2H-1}drdl\\
&=2H(1-K)2^{-K}\int_s^tl^{2H-1}dl\int_{s'}^{t'}
\left(r^{2H}+l^{2H}\right)^{K-2}r^{2H-1}dr\\
&=-2^{-K}\int_{s}^{t}l^{2H-1}\left\{G_l(t')-G_l(s')\right\}dl\\
&=-2H(K-1)2^{-K}(t'-s')\int_{s}^{t}l^{2H-1}
\left(l^{2H}+\xi_l^{2H}\right)^{K-2}\xi_l^{2H-1}dl.
\end{align*}
Notice that
\begin{align*}
\frac1{\left(l^{2H}+\xi_l^{2H}\right)^{2-K}}&\leq \frac1{l^{2H\alpha(2-K)}\xi_l^{2H(1-\alpha)(2-K)}}\\
&=\frac1{l^{\alpha(4H-1)}\xi_l^{(1-\alpha)(4H-1)}} =\frac1{l^{2H}\xi_l^{2H-1}}
\end{align*}
with $1-\alpha=\frac{2H-1}{4H-1}$ by Young's inequality. We get
\begin{align*}
\int_{s}^{t}l^{2H-1}
\left(l^{2H}+\xi_l^{2H}\right)^{K-2}\xi_l^{2H-1}dl&\leq \int_{s}^{t}
l^{-1}dl\leq \frac{t-s}{s},
\end{align*}
which deduces
\begin{align*}
|E(B_{t}-B_{s})(B_{t'}-B_{s'})|\leq \frac{(t-s)(t'-s')}{s}.
\end{align*}
This completes the proof.
\end{proof}
From the proof of the above lemma we also have
\begin{equation}\label{sec2-1-eq2.9=00}
0\leq -E(B_{t} -B_{s})(B_{t'}-B_{s'})\leq
\frac{(t-s)(t'-s')^{1-H}}{(1-H)s^{1-H}}
\end{equation}
holds for all $t>s\geq t'>s'>0$. In fact, under the notations of proof of the above lemma we have
\begin{align*}
0\leq -\left[G_r(t)-G_r(s)\right]&=2H(1-K)(t-s)\frac{\xi^{2H-1}}{ \left(r^{2H}+\xi^{2H}\right)^{2-K}}\\
&=2H(1-K)(t-s)\frac{\xi^{2H-1}}{ \left(r^{2H}+\xi^{2H}\right)^{1-K}\left(r^{2H}+\xi^{2H}\right)}\\
&\leq (t-s)\frac{\xi^{2H-1}}{r^{2H(1-K)}r^H\xi^H}
\leq \frac{t-s}{r^{3H-1}s^{1-H}}
\end{align*}
by Cauchy's inequality. It follows that
\begin{align*}
\left|E(B_{t}-B_{s})(B_{t'}-B_{s'})\right|
&=2H(1-K)\int_{s'}^{t'}\int_s^t
\left(r^{2H}+l^{2H}\right)^{K-2}r^{2H-1}l^{2H-1}drdl\\
&\leq \int_{s'}^{t'}r^{2H-1}\left|G_r(t)-G_r(s)\right|dr\leq \int_{s'}^{t'}\frac{t-s}{r^{H}s^{1-H}}dr\\ &=\frac1{(1-H)s^{1-H}}(t-s) \left(t'^{1-H}-s'^{1-H}\right)\\
&\leq \frac{(t-s)(t'-s')^{1-H}}{(1-H)s^{1-H}},
\end{align*}
which deduces the estimate~\eqref{sec2-1-eq2.9=00}.

\begin{lemma}\label{lem2.4}
For $2HK=1$ we have
\begin{align}\label{sec3-eq3.001}
&\left|E\left[B_s(B_{t}-B_{s})\right]\right| \leq
C_{H,K}\frac{s}{t}(t-s),\\
\label{sec3-eq3.002}
&\left|E\left[B_r(B_{t}-B_{s})\right]\right|
\leq\frac{r}{s}(t-s),\\
\label{sec3-eq3.003}
&\left|E\left[B_s(B_{t}-B_{r})\right]\right|
\leq 4(t-r),\\ \label{sec3-eq3.004}
&\left|E\left[B_t(B_{s}-B_{r})\right]\right| \leq 2(s-r)
\end{align}
for all $t>s>r>0$.
\end{lemma}
\begin{proof}
Keeping the notation in the proof of Lemma~\ref{lem2.3}. For the estimate~\eqref{sec3-eq3.001} we have
\begin{align*}
\left|E\left[B_s(B_t-B_s)\right]\right|
&=\frac1{2^K}\left|(t^{2H}+s^{2H})^K-(t-s)-2^Ks\right|\\
&=\frac{t}{2^K}\left(1+(2^K-1)x-(1+x^{2H})^K\right)\\
&\leq C_{H,K}tx(1-x)=C_{H,K}\frac{s}{t}(t-s)
\end{align*}
with $x=\frac{s}{t}$ by the identities~\eqref{sec3-eq000-1}
and~\eqref{sec3-eq000-2}.

In order to prove the other estimates we define the function
$g_r :{\mathbb R}_{+}\to {\mathbb R}$ for $r>0$ by
$$
x\mapsto g_r(x)=(r^{2H}+x^{2H})^K.
$$
We then have by mean value theorem,
\begin{align*}
g_r(t)-g_r(s)&=(t-s)\xi_r^{2H-1}(r^{2H}+\xi_r^{2H})^{K-1} =(t-s)\xi_r^{2H-1}G_{\xi_r}(r)
\end{align*}
for some $\xi_r\in (s,t)$, and
\begin{align*}
|\xi)_r^{2H-1}G_{\xi_r}(r)-1|&= 1- \frac{\xi_r^{2H-1}}{(r^{2H}+\xi_r^{2H})^{1-K}}=1-\left( \frac{\xi_r^{2H}}{r^{2H}+\xi_r^{2H}}\right)^{1-K}\\
&\leq 1-\frac{\xi_r^{2H}}{r^{2H}+\xi_r^{2H}}\leq \frac{r^{2H}}{r^{2H}+\xi_r^{2H}}\leq \frac{r^{2H}}{s^{2H}},
\end{align*}
which deduces
\begin{align*}
\left|E\left[B_r(B_{t}-B_{s})\right]\right|
&=\frac1{2^K}\left|g_r(t)-g_r(s)-(t-s)\right|\\
&\leq (t-s)|\xi_r^{2H-1}G_{\xi_r}(r)-1|\leq \frac{r^{2H}}{s^{2H}}(t-s)\leq \frac{r}{s}(t-s).
\end{align*}
This gives the estimate~\eqref{sec3-eq3.002}.

For~\eqref{sec3-eq3.003}, by mean value theorem we have
\begin{align*}
\left|E\left[B_s(B_{t}-B_{r})\right]\right|
&=\frac1{2^K}\left|g_s(t)-g_s(r)-(t-r)+2(s-r)\right|\\
&=\frac1{2^K}\left|(t-r)\xi_s^{2H-1}(s^{2H}+\xi_s^{2H})^{K-1}
-(t-r)+2(s-r)\right|\\
&\leq \frac1{2^K}\left((t-r)\xi_s^{2H-1}(s^{2H}+\xi_s^{2H})^{K-1}
+(t-r)+2(s-r)\right)\\
&\leq \frac1{2^K}(t-r)\left(\xi_s^{2H-1}(s^{2H}+\xi_s^{2H})^{K-1}
+3\right)\\
&\leq \frac4{2^K}(t-r)
\end{align*}
for some $\xi_s\in (r,t)$. Similarly, we also have
\begin{align*}
\left|E\left[B_t(B_s-B_r)\right]\right|
&=\frac1{2^K}\left|g_t(s)-g_t(r)+(s-r)\right|\leq \frac2{2^K}(s-r),
\end{align*}
which obtains~\eqref{sec3-eq3.004}. Thus, we complete the proof.
\end{proof}

Let $\varphi(x,y)$ be the density function of $(B_s,B_r)$ ($s>r>0$).
That is
$$
\varphi(x,y)=\frac1{2\pi\rho}\exp\left\{-\frac{1}{2\rho^2}\left(
rx^2-2\mu xy+sy^2\right)\right\},
$$
where $\mu_{s,r}=E(B_sB_r)$ and $\rho^2_{s,r}=rs-\mu^2$.
\begin{lemma}\label{lem2.6}
Let $f\in C^1({\mathbb R})$ admit compact support. Then we have
\begin{align}\label{eq3.7==00}
|E[f''(B_{s})&f(B_{r})]|\leq
C_{H,K}\frac{s^{1/4}}{r^{1/4}(s-r)}\int_{\mathbb{R}}
f^2(x)\varphi_{s}(x)dx
\end{align}
and
\begin{align}\notag
\bigl|E[f'(B_{s})&f'(B_{r})]+E[f''(B_{s})f(B_{r})]\bigr|\\
\label{eq3.7}
&\leq C_{H,K}\left(\frac{1}{s^{3/4}r^{1/4}}+\frac{1}{ r^{3/4}s^{1/4}\sqrt{s-r}}\right)
\int_{\mathbb{R}} f^2(x)(s-r+x^2)\varphi_{s}(x)dx
\end{align}
for all $s>r>0$ and $2HK=1$, where
$\varphi_{s}(x)=\frac1{\sqrt{2\pi s}}e^{-\frac{x^2}{2s}}$.
\end{lemma}
\begin{proof}
Elementary calculus can show that
\begin{align*}
E[f''(B_{s})f(B_{r})]&=\int_{\mathbb{R}^2}
f(x)f(y)\frac{\partial^{2}}{\partial x^2}
\varphi(x,y)dxdy\\
&=\int_{\mathbb{R}^2} f(x)f(y)\left\{\frac1{\rho^4} (rx-\mu_{s,r}
y)^2-\frac{r}{\rho^2_{s,r}}\right\}\varphi(x,y)dxdy
\end{align*}
and
\begin{align*}
\int_{\mathbb{R}^2}|f(y)|^2& \left|\frac{r}{\rho^2_{s,r}}(x-\frac{\mu_{s,r}}ry)^2 -1\right|^2\varphi(x,y)dxdy\\
&=\int_{\mathbb{R}}|f(y)|^2\varphi_r(y)dy \int_{\mathbb{R}}\left|\frac{r}{\rho^2_{s,r}}(x-\frac{\mu_{s,r}}ry)^2 -1\right|^2\varphi_{\frac{\rho^2_{s,r}}{r}}(x-\frac{\mu_{s,r}}ry)dx\\
&=\int_{\mathbb{R}}|f(y)|^2\varphi_r(y)dy \int_{\mathbb{R}}(u^2 -1)^2\varphi_{1}(u)du\\
&=2\int_{\mathbb{R}}|f(y)|^2\varphi_r(y)dy.
\end{align*}
We have
\begin{align*}
|E[&f''(B_{s})f(B_{r})]|\leq \frac{r}{\rho^2_{s,r}}\int_{\mathbb{R}^2} \left|f(x)f(y)\left\{\frac{r}{\rho^2_{s,r}}(x-\frac{\mu_{s,r}}ry)^2 -1\right\}\right|\varphi(x,y)dxdy\\
&\leq \frac{r}{\rho^2_{s,r}}\left(\int_{\mathbb{R}^2} |f(x)|^2\varphi(x,y)dxdy
\int_{\mathbb{R}^2}|f(y)|^2 \left|\frac{r}{\rho^2_{s,r}}(x-\frac{\mu_{s,r}}ry)^2 -1\right|^2\varphi(x,y)dxdy\right)^{1/2}\\
&=C_{H,K}\frac{r}{\rho^2_{s,r}}\left(\int_{\mathbb{R}} |f(x)|^2\varphi_s(x)dx \int_{\mathbb{R}}|f(y)|^2\varphi_r(y)dy\right)^{1/2}\\
&\leq C_{H,K}\frac{r}{\rho^2_{s,r}}\sqrt[4]{\frac{s}{r}}\int_{\mathbb{R}} |f(x)|^2\varphi_s(x)dx\leq C_{H,K}\frac{s^{1/4}}{r^{1/4}(s-r)}\int_{\mathbb{R}} |f(x)|^2\varphi_s(x)dx
\end{align*}
by the inequalities~\eqref{sec2-1-eq2.7} and the fact
\begin{align}\label{eq4.5}
E\left[f^2(B_{r})\right]&=\int_{\mathbb R}f(x)^2\varphi_r(x)dx\\   \tag*{} &\leq
\sqrt{\frac{s}{r}}\int_{\mathbb R}f(x)^2\varphi_s(x)dx=\sqrt{\frac{s}{r}}E\left[f^2(B_{s}) \right]
\end{align}
with $s\geq r>0$.

On the other hand, we have
\begin{align*}
E[f'(B_{s})&f'(B_{r})]=\int_{\mathbb{R}^2}
f(x)f(y)\frac{\partial^{2}}{\partial x\partial
y}\varphi(x,y)dxdy\\
&=\int_{\mathbb{R}^2} f(x)f(y)\left\{\frac1{\rho^4_{s,r}}(sy-\mu_{s,r}x)(rx-\mu_{s,r}
y)+\frac{\mu_{s,r}}{\rho^2_{s,r}}\right\}\varphi(x,y)dxdy,
\end{align*}
which deduce, by the following identity
\begin{align} (sy-\mu
x)(rx-\mu_{s,r}y)=\rho^2_{s,r}y(x-\frac{\mu_{s,r}}{r}y)-\mu_{s,r} r(x-\frac{\mu_{s,r}}{r}y)^2,
\end{align}
\begin{align*}
E[f'(B_{s})&f'(B_{r})]+E[f''(B_{s})f(B_{r})]\\
&= \frac1{\rho^4_{s,r}}\int_{\mathbb{R}^2} f(x)f(y)\left[(sy-\mu_{s,r}
x)(rx-\mu_{s,r}y)+(rx-\mu_{s,r}y)^2\right]\varphi(x,y)dxdy\\
&\qquad\qquad+\frac{\mu_{s,r}-r}{\rho^2_{s,r}}\int_{\mathbb{R}^2}
f(x)f(y)\varphi(x,y)dxdy\\
&=\frac1{\rho^2_{s,r}} \int_{\mathbb{R}^2} f(x)f(y)
y(x-\frac{\mu_{s,r}}{r}y)\varphi(x,y)dxdy \\
&\qquad+\frac{(r-\mu_{s,r})r}{\rho^4_{s,r}} \int_{\mathbb{R}^2}
f(x)f(y)(x-\frac{\mu_{s,r}}{r}y)^2\varphi(x,y)dxdy \\
&\qquad\qquad+\frac{r-\mu_{s,r}}{\rho^2_{s,r}}\int_{\mathbb{R}^2}
f(x)f(y)\varphi(x,y)dxdy\equiv \Lambda_1+\Lambda_2+\Lambda_3.
\end{align*}
Notice that
\begin{align*}
\int_{\mathbb{R}^2}
f^2(y)&|x-\frac{\mu_{s,r}}{r}y|^{2m}\varphi(x,y)dxdy\\
&=\int_{\mathbb{R}}f^2(y)\varphi_{r}(y)dy\int_{\mathbb R}|x-\frac{\mu_{s,r}}{r}y|^{2m}
\frac{\sqrt{r}}{\sqrt{2\pi}\rho_{s,r}}e^{-\frac{r}{2\rho_{s,r}^2} (x-\frac{\mu_{s,r}}{r}y)^2}dx\\
&=C_m(\frac{\rho^2_{s,r}}r)^m\int_{\mathbb{R}}f^2(y)\varphi_{r}(y)dy,
\end{align*}
for all $m\geq 0$, and by the inequalities~\eqref{sec2-1-eq2.7}
\begin{align*}
\int_{\mathbb{R}^2}|f(x)y|^2&\varphi(x,y)dxdy
=\int_{\mathbb{R}}
f^2(x)\varphi_{s}(x)dx\int_{\mathbb{R}}
y^2\frac{\sqrt{s}}{\sqrt{2\pi}\rho_{s,r}}e^{-\frac{s}{2\rho^2_{s,r}} (y-\frac{\mu_{s,r}}{s}x)^2}dy\\
&=\int_{\mathbb{R}}
f^2(x)\varphi_{s}(x)dx\int_{\mathbb{R}}
\left((y-\frac{\mu_{s,r}}sx)+\frac{\mu_{s,r}}sx\right)^2 \frac{\sqrt{s}}{\sqrt{2\pi}\rho_{s,r}}e^{-\frac{s}{2\rho^2_{s,r}} (y-\frac{\mu_{s,r}}{s}x)^2}dy\\
&=\int_{\mathbb{R}}
f^2(x)\varphi_{s}(x)dx\int_{\mathbb{R}}
\left((y-\frac{\mu_{s,r}}sx)^2+\frac{\mu^2_{s,r}}{s^2}x^2\right) \frac{\sqrt{s}}{\sqrt{2\pi}\rho_{s,r}}e^{-\frac{s}{2\rho^2_{s,r}} (y-\frac{\mu_{s,r}}{s}x)^2}dy\\
&=\int_{\mathbb{R}}
f^2(x)\left(\frac{\rho^2_{s,r}}s+\frac{\mu^2_{s,r}}{s^2}x^2\right) \varphi_{s}(x)dx\\
&\leq C_{H,K}\frac{r}{s}\int_{\mathbb{R}}
f^2(x)\left(s-r+x^2\right) \varphi_{s}(x)dx
\end{align*}
for all $0<r<s$. We see that, by the fact~\eqref{eq4.5}
\begin{align*}
|\Lambda_1|&\leq \frac{1}{\rho^2_{s,r}}\left(\int_{\mathbb{R}^2}
f^2(y)(x-\frac{\mu_{s,r}}{r}y)^2 \varphi(x,y)dxdy \int_{\mathbb{R}^2}f^2(x)y^2
\varphi(x,y)dxdy\right)^{1/2}\\
&\leq  \frac{C_{H,K}}{\rho^2_{s,r}}\left(\frac{\rho^2_{s,r}}r\int_{\mathbb{R}}
f^2(y)\varphi_{r}(y)dy\cdot\frac{r}{s}\int_{\mathbb{R}}
f^2(x)(s-r+x^2)\varphi_{s}(x)dx\right)^{1/2}\\
&\leq \frac{ C_{H,K}}{\sqrt{s}\rho_{s,r}}\left(\sqrt{\frac{s}{r}}\int_{\mathbb{R}}
f^2(y)y^2\varphi_{s}(y)dy\int_{\mathbb{R}}
f^2(x)(s-r+x^2)\varphi_{s}(x)dx\right)^{1/2}\\
&\leq \frac{C_{H,K}}{(rs)^{1/4}\rho_{s,r}}\int_{\mathbb{R}}
f^2(x)(s-r+x^2)\varphi_{s}(x)dx,
\end{align*}
\begin{align*}
|\Lambda_2|&\leq \frac{(r-\mu_{s,r})r}{\rho^4_{s,r}} \left(\int_{\mathbb{R}^2}
f^2(y)(x-\frac{\mu_{s,r}}{r}y)^4 \varphi(x,y)dxdy
\int_{\mathbb{R}^2}f^2(x)
\varphi(x,y)dxdy\right)^{1/2}\\
&=\frac{\sqrt{3}(r-\mu_{s,r})}{\rho^2_{s,r}}\sqrt[4]{\frac{s}{r}}
\int_{\mathbb{R}}f^2(x) \varphi_{s}(x)dx
\end{align*}
and
\begin{align*}
|\Lambda_3|&\leq \frac{r-\mu_{s,r}}{\rho^2_{s,r}} \left(\int_{\mathbb{R}^2}
f^2(y) \varphi(x,y)dxdy
\int_{\mathbb{R}^2}f^2(x)
\varphi(x,y)dxdy\right)^{1/2}\\
&=\frac{r-\mu_{s,r}}{\rho^2_{s,r}}\sqrt[4]{\frac{s}{r}}
\int_{\mathbb{R}}f^2(x) \varphi_{s}(x)dx.
\end{align*}
Thus, the estimate~\eqref{eq3.7} follows from Lemma~\ref{lem2.2} and Lemma~\ref{lem2.1}.
\end{proof}
From the above proof of Lemma~\ref{lem2.6} we also have
\begin{align}\label{eq3.7-000}
|E[f(B_{r})&f''(B_{s})]|\leq
\frac{C_{H,K}r^{5/4}}{s^{5/4}(r-s)}\int_{\mathbb{R}}
f^2(x)\varphi_r(x)dx
\end{align}
and
\begin{align}\notag
\bigl|E[f'(B_{s})&f'(B_{r})]+E[f(B_{r})f''(B_{s})]\bigr|\\
\label{eq3.7-00} &\leq
C_{H,K}\left(\frac{r^{1/4}}{s^{5/4}\sqrt{r-s}}+\frac{r^{1/4}}{
s^{5/4}}\right)\int_{\mathbb{R}}
f^2(x)\left(r-s+x^2\right)\varphi_{r}(x)dx
\end{align}
for all $0<s<r$ and $2HK=1$.

\section{Existence of quadratic covariation}\label{sec3}
In this section, we study the quadratic covariation
$\left[f(B),B\right]$. Denote
$$
J_\varepsilon(f,t):=\frac{1}{\varepsilon}\int_0^t\left\{ f(B_{
s+\varepsilon}) -f(B_s)\right\}(B_{s+\varepsilon}-B_s)ds
$$
for $\varepsilon>0$ and $0\leq t\leq T$. Recall that the quadratic
covariation, the forward integral and the backward integrals are
defined as
\begin{align}\label{sec4-eq1.1}
[f(B),B]_t:&=\lim_{\varepsilon\downarrow
0}J_\varepsilon(f,t),\\
\int_0^tf(B_s)d^{-}B_s:&=\lim_{\varepsilon\to 0}
\frac{1}{\varepsilon}\int_0^tf(B_s)
(B_{s+\varepsilon}-B_s)ds,\\
\int_0^tf(B_s)d^{+}B_s:&=\lim_{\varepsilon\to 0}
\frac{1}{\varepsilon}\int_0^tf(B_{
s+\varepsilon})(B_{s+\varepsilon}-B_s)ds,
\end{align}
provided the corresponding limits exist in $L^1$, and we have
\begin{align}
[f(B),B]_t=\int_0^tf'(B_s)d[B,B]_s=2^{1-K}\int_0^tf'(B_s)ds
\end{align}
for all $0\leq t\leq T$ and $f\in C^1({\mathbb R})$ (see
Russo-Tudor~\cite{Russo-Tudor} and
Russo-Vallois~\cite{Russo-Vallois2,Russo-Vallois3}).

Now, we study the existence in $L^2$ of the forward integral,
backward integral and quadratic covariation. Consider the set
${\mathscr H}$ of measurable functions $f$ on $\mathbb R$ such that
$\|f\|_{\mathscr H}<+ \infty$, where
\begin{align*}
\|f\|_{\mathscr H}^2:=\int_0^T\int_{\mathbb
R}|f(x)|^2\varphi_s(x)dxds+\int_0^T\int_{\mathbb
R}|f(x)x|^2\varphi_s(x)\frac{dxds}{s}
\end{align*}
with $\varphi_s(x)=\frac{1}{\sqrt{2\pi s}}e^{-\frac{x^2}{2s}}$.
Clearly, ${\mathscr H}$ is a Banach space and the set ${\mathscr E}$
of elementary functions
$$
f_\Delta(x)=\sum_jf_j1_{(a_j,a_{j+1}]}(x),\quad f_j\in {\mathbb
R};-\infty<a_0<a_1<\cdots<a_N<\infty
$$
is dense in ${\mathscr H}$, and moreover every $f\in {\mathscr H}$
is locally square integrable and the space of  measurable functions
\begin{align*}
{\mathscr H}_K=\left\{f\;\;|\;\;\int_0^T\int_{\mathbb
R}|f(x)|^{{2}/{K}}\varphi_s(x)dxds<\infty\right\}
\end{align*}
is a simple subspace of ${\mathscr H}$.
\begin{lemma}\label{lemma4.1-0}
Let $2HK=1$. If $f\in {\mathscr H}$, then we have
\begin{align}\label{sec4-eq3-2}
E\left|\frac{1}{\varepsilon}\int_0^tf(B_s)
(B_{s+\varepsilon}-B_s)ds\right|^2
\leq C_{H,K}\|f\|^2_{\mathscr H},\\ \label{sec4-eq3-2-1}
E\left|\frac{1}{\varepsilon}\int_0^tf(B_{
s+\varepsilon})(B_{s+\varepsilon}-B_s)ds\right|^2\leq
C_{H,K}\|f\|^2_{\mathscr H}
\end{align}
for all $0<\varepsilon<T$ and $0\leq t\leq T$.
\end{lemma}
\begin{proof}
Without loss of generality one may assume that $T=1$. We prove only
the estimate~\eqref{sec4-eq3-2} and similarly one can
prove~\eqref{sec4-eq3-2-1}. Let $0<\varepsilon<T$ and $0<s,r<T$. By
approximating we may assume that $f$ is an infinitely differentiable
function with compact support. It follows that
\begin{align}\notag
E&\left[f(B_{
s})f(B_{r})(B_{s+\varepsilon}-B_s)(B_{r+\varepsilon}-B_r) \right]\\
\notag &=E\left[f(B_{ s})f(B_{
r})(B_{s+\varepsilon}-B_s)\int_r^{r+\varepsilon}dB_l\right]\\ \notag
&=E\int_0^T\int_0^T1_{[r,r+\varepsilon]}(u)D_v^{H,K}\left[
f(B_{s})f(B_{r})(B_{s+\varepsilon}-B_s)\right]\phi(u,v)dudv \\
\notag &=\left(\int_0^T\int_0^T1_{[r,r+\varepsilon]}(u)1_{[0,s]}(u)
\phi(u,v)dudv\right)E\left[
f'(B_{s})f(B_{r})(B_{s+\varepsilon}-B_s)\right]\\  \notag
&\qquad+\left(\int_0^T\int_0^T1_{[r,r+\varepsilon]}(u) 1_{[0,r]}(u)
\phi(u,v)dudv\right) E\left[
f(B_{s})f'(B_{r})(B_{s+\varepsilon}-B_s)\right]\\  \notag
&\qquad\qquad+\left(\int_0^T\int_0^T1_{[r,r+\varepsilon]}(u)
1_{[s,s+\varepsilon]}(u) \phi(u,v)dudv\right)E\left[
f(B_{s})f(B_{r})\right]\\ \notag
&=E\left[B_{s}(B_{r+\varepsilon}-B_r)\right]E\left[
f'(B_{s})f(B_{r})(B_{s+\varepsilon}-B_s) \right]\\
\label{lem4.1-eq1} &\qquad
+E\left[B_{r}(B_{r+\varepsilon}-B_r)\right]E\left[
f(B_{s})f'(B_{r})(B_{s+\varepsilon}-B_s) \right]\\ \notag
&\qquad\qquad +E\left[(B_{r+\varepsilon}-B_r)
(B_{s+\varepsilon}-B_s)\right]E\left[ f(B_{s})f(B_{r})\right]\\
\notag &\equiv \Psi_{\varepsilon}(s,r,1)+\Psi_{\varepsilon}(s,r,2)
+\Psi_{\varepsilon}(s,r,3).
\end{align}
In order to establish~\eqref{sec4-eq3-2} we first show that
\begin{equation}\label{eq4.4}
\frac{1}{\varepsilon^2}\left|\int_0^t\int_0^t\Psi_{\varepsilon}
(s,r,3)ds dr \right|\leq C_{H,K}\|f\|^2_{\mathscr H}
\end{equation}
for all $\varepsilon>0$ small enough. we have
\begin{align*}
\frac{1}{\varepsilon^{2}}&\left|\int_0^1\int_0^1
\Psi_{\varepsilon}(s,r,3) ds dr \right|\\
&\leq \frac{2}{\varepsilon^{2}}\int_0^1\int_0^s
|E\left[(B_{r+\varepsilon}-B_r) (B_{s+\varepsilon}-B_s)\right]|
E\left|f(B_{s}) f(B_{r})\right|dsdr\\
&=\frac{2}{\varepsilon^{2}}\int_\varepsilon^1ds \int_0^{s-\varepsilon}|E\left[(B_{r+\varepsilon}-B_r) (B_{s+\varepsilon}-B_s)\right]|
E\left|f(B_{s})f(B_{r})\right|dr\\
&\qquad+\frac{2}{\varepsilon^{2}}\int_0^\varepsilon\int_0^s
|E\left[(B_{r+\varepsilon}-B_r) (B_{s+\varepsilon}-B_s)\right]|
E\left|f(B_{s}) f(B_{r})\right|dsdr\\
&\qquad+\frac{2}{\varepsilon^{2}}\int_\varepsilon^1\int_{s-\varepsilon}^s
|E\left[(B_{r+\varepsilon}-B_r) (B_{s+\varepsilon}-B_s)\right]|
E\left|f(B_{s}) f(B_{r})\right|dsdr\\
&\equiv \Lambda_{31}+\Lambda_{32}+\Lambda_{33}
\end{align*}
for all $0<\varepsilon\leq 1$. Clearly, Lemma~\eqref{lem2.3} and the
fact~\eqref{eq4.5} imply that
\begin{align*}
\Lambda_{31}&\leq 2\int_\varepsilon^1ds \int_0^{s-\varepsilon}\frac1s
\sqrt{E(|f(B_{s})|^2)E(|f(B_{r})|^2)}dr\\
&\leq 2\int_0^1ds \int_0^{s}\frac1{s^{3/4}r^{1/4}}E(|f(B_{s})|^2)dr\\
&=2\int_0^1dsE(|f(B_{s})|^2)dr\leq 2\|f\|_{\mathscr H}^2.
\end{align*}
Notice that
\begin{align*}
\int_0^\varepsilon ds\int_0^s\left|Ef(B_{s}) f(B_{r})\right|dr &\leq \int_0^\varepsilon ds\int_0^s\sqrt{E(|f(B_{s})|^2)E(|f(B_{r})|^2)}dr\\
&\leq \int_0^\varepsilon ds\int_0^s\frac{s^{1/4}}{r^{1/4}}E(|f(B_{s})|^2)dr\\
&\leq \frac43\varepsilon\int_0^\varepsilon E(|f(B_{s})|^2)ds\leq \varepsilon\|f\|_{\mathscr H}^2
\end{align*}
and
\begin{align*}
\int_\varepsilon^1ds\int_{s-\varepsilon}^s\left|Ef(B_{s}) f(B_{r})\right|dr
&\leq \int_\varepsilon^1ds\int_{s-\varepsilon}^s
\sqrt{E(|f(B_{s})|^2)E(|f(B_{r})|^2)}dr\\
&\leq \int_\varepsilon^1ds\int_{s-\varepsilon}^s \frac{s^{1/4}}{r^{1/4}}E(|f(B_{s})|^2)dr\\
&=\frac43\int_\varepsilon^1E(|f(B_{s})|^2) s^{1/4}\left(s^{3/4}-(s-\varepsilon)^{3/4}\right)ds \\
&\leq \frac43\varepsilon\int_\varepsilon^1E(|f(B_{s})|^2)ds\leq \frac43\varepsilon\|f\|_{\mathscr H}^2
\end{align*}
for all $0<\varepsilon\leq 1$. We get
\begin{align*}
\Lambda_{32}+\Lambda_{33}&\leq
\frac{2}{\varepsilon^{2}}\int_0^\varepsilon\int_0^s
\sqrt{E[(B_{r+\varepsilon}-B_r)^2]E[(B_{s+\varepsilon}-B_s)^2]}
E\left|f(B_{s}) f(B_{r})\right|dsdr\\
&\qquad+\frac{2}{\varepsilon^{2}} \int_\varepsilon^1\int_{s-\varepsilon}^s
\sqrt{E[(B_{r+\varepsilon}-B_r)^2]E[(B_{s+\varepsilon}-B_s)^2]}
E\left|f(B_{s})f(B_{r})\right|dsdr\\
&\leq C_{H,K}\frac{1}{\varepsilon}\int_0^\varepsilon\int_0^sE\left|f(B_{s}) f(B_{r})\right|dsdr+C_{H,K}\frac{1}{\varepsilon} \int_\varepsilon^1\int_{s-\varepsilon}^s
E\left|f(B_{s})f(B_{r})\right|dsdr\\
&\leq C_{H,K}\|f\|_{\mathscr H}^2.
\end{align*}
It follows that
\begin{align*}
\frac{1}{\varepsilon^{2}}&\left|\int_0^1\int_0^1
\Psi_{\varepsilon}(s,r,3) ds dr \right|\leq \Lambda_{31}+\Lambda_{32}+\Lambda_{33}\leq C_{H,K}\|f\|_{\mathscr H}^2
\end{align*}
for all $0<\varepsilon\leq 1$.

Now, let us prove
\begin{equation}\label{eq4.4-00}
\frac{1}{\varepsilon^2}\int_0^t\int_0^t\left|\Psi_{\varepsilon}
(s,r,1)+\Psi_{\varepsilon} (s,r,2)\right| drds \leq
C_{H,K}\|f\|^2_{\mathscr H}
\end{equation}
for all $\varepsilon>0$. We have
\begin{align}\notag
\Psi_{\varepsilon}(s,r,1)
&=E\left[B_{s}(B_{r+\varepsilon}-B_r)\right]E\left[
f'(B_{s})f(B_{r})(B_{s+\varepsilon}-B_s) \right]\\ \label{decom1-11}
&=E\left[B_{s}(B_{r+\varepsilon}-B_r)\right]
E\left[B_{s}(B_{s+\varepsilon}-B_s)\right]
E\left[f''(B_{s})f(B_{r})\right]\\  \label{decom1-21}
&\quad+E\left[B_{s}(B_{r+\varepsilon}-B_r)\right]
E\left[B_{r}(B_{s+\varepsilon}-B_s)\right]
E\left[f'(B_{s})f'(B_{r})\right],\\
\notag \Psi_{\varepsilon}(s,r,2)&=
E\left[B_{r}(B_{r+\varepsilon}-B_r)\right]E\left[
f(B_{s})f'(B_{r})(B_{s+\varepsilon}-B_s) \right]\\
\label{decom2-11} &=E\left[B_{r}(B_{r+\varepsilon}-B_r)\right]
E\left[B_{s}(B_{s+\varepsilon}-B_s)\right]E\left[
f'(B_{s})f'(B_{r})\right]\\    \label{decom2-21}
&\qquad+E\left[B_{r}(B_{r+\varepsilon}-B_r)\right]
E\left[B_{r}(B_{s+\varepsilon}-B_s)\right]E\left[
f(B_{s})f''(B_{r})\right].
\end{align}
For $s>r>0$ we decompose
$\Psi_{\varepsilon}(s,r,1)+\Psi_{\varepsilon}(s,r,2)$ as follows
\begin{align*}
\Psi_{\varepsilon}(s,r,1)
=&E\left[B_{s}(B_{r+\varepsilon}-B_r)\right]
E\left[B_{r}(B_{s+\varepsilon}-B_s)\right]\\
&\cdot\Bigl(E\left[f''(B_{s})f(B_{r})\right]
+E\left[f'(B_{s})f'(B_{r})\right]\Bigr)\\
&\quad+\Bigl\{E\left[B_{s}(B_{r+\varepsilon}-B_r)\right]
E\left[B_{s}(B_{s+\varepsilon}-B_s)\right]\\
&\qquad\qquad- E\left[B_{s}(B_{r+\varepsilon}-B_r)\right]
E\left[B_{r}(B_{s+\varepsilon}-B_s)\right] \Bigr\}
E\left[f''(B_{s})f(B_{r})\right]\\
&\equiv\Psi_{\varepsilon}(s>r,1)+\Psi_{\varepsilon}(s>r,2)
\end{align*}
and
\begin{align*}
\Psi_{\varepsilon}(s,r,2)
&=E\left[B_{r}(B_{r+\varepsilon}-B_r)\right] E\left[B_{s}(B_{s+\varepsilon}-B_s)\right]\\
&\qquad\cdot\Bigl\{E\left[
f'(B_{s})f'(B_{r})\right]+E\left[f(B_{s})f''(B_{r})\right]\Bigr\}\\
&\qquad+\Bigl\{E\left[B_{r}(B_{r+\varepsilon}-B_r)\right]
E\left[B_{r}(B_{s+\varepsilon}-B_s)\right]\\
&\qquad\qquad- E\left[B_{r}(B_{r+\varepsilon}-B_r)\right]
E\left[B_{s}(B_{s+\varepsilon}-B_s)\right]\Bigr\}E\left[
f(B_{s})f''(B_{r})\right]\\
&\equiv\Psi_{\varepsilon}(s>r,3)+\Psi_{\varepsilon}(s>r,4).
\end{align*}
Notice that
\begin{align*}
|\Psi_{\varepsilon}(s>r,1)|&\leq C_{H,K}\frac{r}{s}\varepsilon^2
\left|E[f''(B_{s})f(B_{r})]+E[f'(B_{s})f'(B_{r})]\right|,\\
|\Psi_{\varepsilon}(s>r,2)|&= \left|E[B_{s}(B_{r+\varepsilon}-B_r)]
E[(B_{s}-B_{r})(B_{s+\varepsilon}-B_s)]
E[f''(B_{s})f(B_{r})t]\right|\\
&\leq C_{H,K}\varepsilon^2\frac{s-r}{s}\left|E[f''(B_{s})f(B_{r})
]\right|,\\
|\Psi_{\varepsilon}(s>r,3)|&\leq \varepsilon^2\left|E[
f'(B_{s})f'(B_{r})]+E[f(B_{s})f''(B_{r})]\right|\\
|\Psi_{\varepsilon}(s>r,4)|& =\left|E[B_{r}(B_{r+\varepsilon}-B_r)]
E[(B_{r}-B_{s})(B_{s+\varepsilon}-B_s)]E[
f(B_{s})f''(B_{r})]\right|\\
&\leq C_{H,K}\varepsilon^2\frac{s-r}{s}\left|E[
f(B_{s})f''(B_{r})]\right|
\end{align*}
by Lemma~\ref{lem2.3} and Lemma~\ref{lem2.4}. We get
\begin{align*}
\frac{1}{\varepsilon^2}\int_0^tds\int_0^s &|\Psi_{\varepsilon}
(s,r,1)+\Psi_{\varepsilon} (s,r,2)|dr\\
&\leq \frac{1}{\varepsilon^2}\sum_{i=1}^4\int_0^tds
\int_0^s|\Psi_{\varepsilon}(s>r,i)|dr\leq C_{H,K}\|f\|^2_{\mathscr
H}
\end{align*}
by Lemma~\ref{lem2.6}. Similarly, for $r>s>0$ in order to decompose
$\Psi_{\varepsilon}(s,r,1)+\Psi_{\varepsilon}(s,r,2)$, we have
\begin{align*}
\eqref{decom1-21}+\eqref{decom2-21}
&=E\left[B_{s}(B_{r+\varepsilon}-B_r)\right]
E\left[B_{r}(B_{s+\varepsilon}-B_s)\right]\\
&\qquad\qquad\cdot\left\{E\left[f'(B_{s})f'(B_{r})\right]+E\left[
f(B_{s})f''(B_{r})\right]\right\}\\
&\quad+\Bigl\{E\left[B_{r}(B_{r+\varepsilon}-B_r)\right]
E\left[B_{r}(B_{s+\varepsilon}-B_s)\right]\\
&\qquad\qquad-E\left[B_{s}(B_{r+\varepsilon}-B_r)\right]
E\left[B_{r}(B_{s+\varepsilon}-B_s)\right]\Bigr\}E\left[
f(B_{s})f''(B_{r})\right]\\
&\equiv \Psi_\varepsilon(s<r,1)+\Psi_\varepsilon(s<r,2)
\end{align*}
and
\begin{align*}
\eqref{decom1-11}+~\eqref{decom2-11}
&=E\left[B_{r}(B_{r+\varepsilon}-B_r)\right]
E\left[B_{s}(B_{s+\varepsilon}-B_s)\right]\\
&\qquad\qquad\cdot\Bigl\{E\left[
f'(B_{s})f'(B_{r})\right]+E\left[f''(B_{s})f(B_{r})\right]\Bigr\}\\
&\quad+\Bigl\{E\left[B_{s}(B_{r+\varepsilon}-B_r)\right]
E\left[B_{s}(B_{s+\varepsilon}-B_s)\right]\\
&\qquad\qquad-E\left[B_{r}(B_{r+\varepsilon}-B_r)\right]
E\left[B_{s}(B_{s+\varepsilon}-B_s)\right]\Bigr\}
E\left[f''(B_{s})f(B_{r})\right]\\
&\equiv \Psi_\varepsilon(s<r,3)+\Psi_\varepsilon(s<r,4),
\end{align*}
which gives
$$
\Psi_{\varepsilon}(s,r,1)+\Psi_{\varepsilon}(s,r,2)
=\sum_{i=1}^4\Psi_{\varepsilon}(s<r,i)
$$
Clearly, Lemma~\ref{lem2.3} and Lemma~\ref{lem2.4} implies that
\begin{align*}
|\Psi_\varepsilon(s<r,1)|&\leq C_{H,K}\varepsilon^2\frac{s}{r}
\left|E[f'(B_{s})f'(B_{r})]+E[f(B_{s})f''(B_{r})]\right|,\\
|\Psi_\varepsilon(s<r,2)|&=\left|E[B_{r}(B_{s+\varepsilon}-B_s)]
E[(B_{r}-B_{s})(B_{r+\varepsilon}-B_r)]E[f(B_{s})f''(B_{r})]\right|\\
&\leq \varepsilon^2\frac{r-s}r\left|E[ f(B_{s})f''(B_{r})]\right|\\
|\Psi_\varepsilon(s<r,3)|&\leq C_{H,K}\varepsilon^2 \left|E[
f'(B_{s})f'(B_{r})]+E[f''(B_{s})f(B_{r})]\right|\\
\Psi_\varepsilon(s<r,4)&\leq \left|E[B_{s}(B_{s+\varepsilon}-B_s)]
E[(B_{r}-B_{s})(B_{r+\varepsilon}-B_r)]
E[f''(B_{s})f(B_{r})]\right|\\
&\leq \varepsilon^2\frac{r-s}{s}\left|E[f''(B_{s})f(B_{r})]\right|
\end{align*}
for all $r>s>0$. It follows from~\eqref{eq3.7-000}
and~\eqref{eq3.7-00} that
\begin{align*}
\frac{1}{\varepsilon^2}\int_0^tdr\int_0^r &|\Psi_{\varepsilon}
(s,r,1)+\Psi_{\varepsilon} (s,r,2)|ds\\
&\leq \frac{1}{\varepsilon^2}\sum_{i=1}^4\int_0^tds
\int_0^s|\Psi_{\varepsilon}(s<r,i)|dr\leq C_{H,K}\|f\|^2_{\mathscr
H}.
\end{align*}
Thus, we have given the desired estimate~\eqref{eq4.4-00}, and the
lemma follows.
\end{proof}

In this section our main result is the following theorem which shows
that $J_{\varepsilon}(f,t)$ converges in $L^2$ as $\varepsilon$
tends to $0$.
\begin{theorem}\label{th3.1}
Let $2HK=1$. If $f\in {\mathscr H}$, then the forward, backward
integrals $\int_0^tf(B_s)d^{\mp}B_s$ and the quadratic covariation
$[f(B),B]$ exist in $L^2$, and
\begin{align}
E\left|\int_0^tf(B_s)d^{\pm}B_s\right|^2\leq
C_{H,K}\|f\|^2_{\mathscr H}
\end{align}
for all $0\leq t\leq T$.
\end{theorem}
\begin{proof}
From Lemma~\ref{lemma4.1-0}, it is enough to show that
\begin{equation}\label{sec40-eq3-1}
E\left|J_{\varepsilon_1}^{-}-J_{\varepsilon_2}^{-}\right|^2
\longrightarrow 0,
\end{equation}
and
\begin{equation}\label{sec40-eq3-2}
E\left|J_{\varepsilon_1}^{+}-J_{\varepsilon_2}^{+}\right|^2
\longrightarrow 0
\end{equation}
as $\varepsilon_1,\varepsilon_2\downarrow 0$, where
\begin{align*}
J_{\varepsilon}^{-}=\frac{1}{\varepsilon}\int_0^tf(B_s)
(B_{s+\varepsilon}-B_s)ds \quad {\text {and }}\quad
J_{\varepsilon}^{+}=\frac{1}{\varepsilon}\int_0^tf(B_{
s+\varepsilon})(B_{s+\varepsilon}-B_s)ds.
\end{align*}
Without loss of generality we assume that $\varepsilon_1>\varepsilon_2$. We prove only the convergence~\eqref{sec40-eq3-1} and similarly one
can prove~\eqref{sec40-eq3-2}. It follows that
\begin{align*}
E\bigl|J_{\varepsilon_1}^{-}&-J_{\varepsilon_2}^{-}\bigr|^2
=\frac1{\varepsilon_1^2} \int_0^t\int_0^tEf(B_s)f(B_r)
(B_{s+\varepsilon_1}-B_s) (B_{r+\varepsilon_1}-B_r)dsdr\\
&\qquad-2
\frac1{\varepsilon_1\varepsilon_2} \int_0^t\int_0^tEf(B_s)f(B_r)
(B_{s+\varepsilon_1}-B_s) (B_{r+\varepsilon_2}-B_r)dsdr\\
&\qquad+\frac1{\varepsilon_2^2} \int_0^t\int_0^tEf(B_s)f(B_r)
(B_{s+\varepsilon_2}-B_s) (B_{r+\varepsilon_2}-B_r)dsdr\\
&\equiv \frac1{\varepsilon_1^2\varepsilon_2}\int_0^t\int_0^t
\left\{\varepsilon_2\Phi_{s,r}(1,\varepsilon_1)-\varepsilon_1
\Phi_{s,r}(2,\varepsilon_1,\varepsilon_2)\right\}dsdr\\
&\qquad+
\frac1{\varepsilon_1\varepsilon_2^2}\int_0^t\int_0^t \left\{
\varepsilon_1\Phi_{s,r}(1,\varepsilon_2)
-\varepsilon_2
\Phi_{s,r}(2,\varepsilon_1,\varepsilon_2)\right\}dsdr,
\end{align*}
where
$$
\Phi_{s,r}(1,\varepsilon) =E\left[f(B_s)f(B_r)(B_{s+\varepsilon}
-B_s) (B_{r+\varepsilon}-B_r)\right],
$$
and
$$
\Phi_{s,r}(2,\varepsilon_1,\varepsilon_2) =E\left[f(B_s)f(B_r)
(B_{s+\varepsilon_1}-B_s) (B_{r+\varepsilon_2}-B_r)\right].
$$
We have by~\eqref{lem4.1-eq1}
\begin{align*}
\Phi_{s,r}(1,\varepsilon) &=\Psi_{\varepsilon}(s,r,1)
+\Psi_{\varepsilon}(s,r,2)
+\Psi_{\varepsilon}(s,r,3)\\
&=E\left[B_{s}(B_{r+\varepsilon}-B_r)\right]
E\left[B_{s}(B_{s+\varepsilon}-B_s)\right]
E\left[f''(B_{s})f(B_{r})\right]\\
&+E\left[B_{s}(B_{r+\varepsilon}-B_r)\right]
E\left[B_{r}(B_{s+\varepsilon}-B_s)\right]
E\left[f'(B_{s})f'(B_{r})\right]\\
&\qquad+E\left[B_{r}(B_{r+\varepsilon}-B_r)\right]
E\left[B_{s}(B_{s+\varepsilon}-B_s)\right]
E\left[f'(B_{s})f'(B_{r})\right]\\
&\qquad+E\left[B_{r}(B_{r+\varepsilon}-B_r)\right]
E\left[B_{r}(B_{s+\varepsilon}-B_s)\right]
E\left[f(B_{s})f''(B_{r})\right]\\
&\qquad\qquad+E\left[(B_{r+\varepsilon}-B_r)
(B_{s+\varepsilon}-B_s)\right]E\left[
f(B_{s})f(B_{r})\right]
\end{align*}
and
\begin{align*}
\Phi_{s,r}(2,\varepsilon_1,\varepsilon_2)&=
E\left[B_{s}(B_{r+\varepsilon_2}-B_r)\right]E\left[
f'(B_{s})f(B_{r})(B_{s+\varepsilon_1}-B_s) \right]\\
&\qquad +E\left[B_{r}(B_{r+\varepsilon_2}-B_r)\right]E\left[
f(B_{s})f'(B_{r})(B_{s+\varepsilon_1}-B_s) \right]\\
&\qquad\qquad +E\left[(B_{s+\varepsilon_1}-B_s)(B_{r+\varepsilon_2}-B_r)
\right]E\left[
f(B_{s})f(B_{r})\right]\\
&=E\left[B_{s}(B_{r+\varepsilon_2}-B_r)\right]
E\left[B_{s}(B_{s+\varepsilon_1}-B_s)\right]
E\left[f''(B_{s})f(B_{r})\right]\\
&+E\left[B_{s}(B_{r+\varepsilon_2}-B_r)\right]
E\left[B_{r}(B_{s+\varepsilon_1}-B_s)\right]
E\left[f'(B_{s})f'(B_{r})\right]\\
&\qquad +E\left[B_{r}(B_{r+\varepsilon_2}-B_r)\right]
E\left[B_{s}(B_{s+\varepsilon_1}-B_s)\right]
E\left[f'(B_{s})f'(B_{r})\right]\\
&\qquad +E\left[B_{r}(B_{r+\varepsilon_2}-B_r)\right]
E\left[B_{r}(B_{s+\varepsilon_1}-B_s)\right]
E\left[f(B_{s})f''(B_{r})\right]\\
&\qquad\qquad +E\left[(B_{s+\varepsilon_1}-B_s)(B_{r+\varepsilon_2}-B_r)
\right]E\left[f(B_{s})f(B_{r})\right].
\end{align*}
Denote
\begin{align*}
A_1(s,r,\varepsilon,j):&=\varepsilon_j E\left[B_{s}(B_{r+\varepsilon}-B_r)\right]
E\left[B_{s}(B_{s+\varepsilon}-B_s)\right]\\
&\qquad-\varepsilon E\left[B_{s}(B_{r+\varepsilon_2}-B_r)\right]
E\left[B_{s}(B_{s+\varepsilon_1}-B_s)\right]\\
A_{21}(s,r,\varepsilon,j):&=\varepsilon_j
E\left[B_{s}(B_{r+\varepsilon}-B_r)\right]
E\left[B_{r}(B_{s+\varepsilon}-B_s)\right]\\
&\qquad-\varepsilon
E\left[B_{s}(B_{r+\varepsilon_2}-B_r)\right]
E\left[B_{r}(B_{s+\varepsilon_1}-B_s)\right]\\
A_{22}(s,r,\varepsilon,j):&=\varepsilon_j
E\left[B_{r}(B_{r+\varepsilon}-B_r)\right] E\left[B_{s}
(B_{s+\varepsilon}-B_s)\right]\\
&\qquad-\varepsilon E\left[B_{r}(B_{r+\varepsilon_2}-B_r)\right]
E\left[B_{s}(B_{s+\varepsilon_1}-B_s)\right]\\
A_3(s,r,\varepsilon,j):&=\varepsilon_j E\left[B_{r}
(B_{r+\varepsilon}-B_r)\right]
E\left[B_{r}(B_{s+\varepsilon}-B_s)\right]\\
&\qquad-\varepsilon E\left[B_{r}(B_{r+\varepsilon_2}-B_r)\right]
E\left[B_{r}(B_{s+\varepsilon_1}-B_s)\right]\\
A_4(s,r,\varepsilon,j):&=\varepsilon_j
E\left[(B_{r+\varepsilon}-B_r)(B_{s+\varepsilon}-B_s)\right]
-\varepsilon
E\left[(B_{s+\varepsilon_1}-B_s)(B_{r+\varepsilon_2}-B_r)\right]
\end{align*}
with $j=1,2$. It follows that
\begin{align*}
&\varepsilon_j\Phi_{s,r}(1,\varepsilon_i) -\varepsilon_i
\Phi_{s,r}(2,\varepsilon_1,\varepsilon_2)\\
&=\Bigl(A_1(s,r,\varepsilon_i,j)E\left[f''(B_{s})f(B_{r})\right]
+\left(A_{21}(s,r,\varepsilon_i,j)+A_{22}(s,r,\varepsilon_i,j)
\right)E\left[f'(B_{s})f'(B_{r})\right]\\
&\qquad+A_3(s,r,\varepsilon_i,j)E\left[f(B_{s})f''(B_{r})
\right]\Bigr)+A_4(s,r,\varepsilon_i,j)E\left[f(B_{s})f(B_{r})
\right]\\
&\equiv\Upsilon(s,r,\varepsilon_i,j)
+A_4(s,r,\varepsilon_i,j)E\left[f(B_{s})f(B_{r}) \right]
\end{align*}
with $i,j=1,2$ and $i\neq j$. In order to end the proof we claim
that the following convergence hold:
\begin{align}\label{sec4-Con-eq1}
\frac1{\varepsilon_i^2\varepsilon_j}\int_0^t\int_0^t
\left\{\varepsilon_j\Phi_{s,r}(1,\varepsilon_i)-\varepsilon_i
\Phi_{s,r}(2,\varepsilon_1,\varepsilon_2)\right\}dsdr
\longrightarrow 0 \qquad (i,j=1,2,i\neq j),
\end{align}
as $\varepsilon_1,\varepsilon_2\to 0$. This will be done in three
parts. Keeping the notations in the proof of Lemma~\ref{lem2.4}.

{\bf Part A}. The following convergence hold:
\begin{align}\label{Part-I-eq1}
\frac1{\varepsilon_i^2\varepsilon_j}\int_0^t\int_0^s
\Upsilon(s,r,\varepsilon_i,j)drds \longrightarrow 0 \qquad
(i,j=1,2,i\neq j)
\end{align}
as $\varepsilon_1,\varepsilon_2\to 0$. For $s>r>0$ we decompose
$\Upsilon(s,r,\varepsilon_i,j)$ as follows
\begin{align*}
\Upsilon(s,r,\varepsilon_i,j)
&=A_1(s,r,\varepsilon_i,j)E\left[f''(B_{s})f(B_{r})\right]
+A_3(s,r,\varepsilon_i,j)E\left[f(B_{s})f''(B_{r})\right]\\
&\qquad+\left(A_{21}(s,r,\varepsilon_i,j)+A_{22}(s,r,\varepsilon_i,j)
\right)E\left[f'(B_{s})f'(B_{r})\right]\\
&=A_{21}(s,r,\varepsilon_i,j)\left\{ E\left[f'(B_{s})f'(B_{r})
\right] +E\left[f''(B_{s})f(B_{r})\right]\right\}\\
&\qquad+A_{22}(s,r,\varepsilon_i,j)\left\{ E\left[f'(B_{s})f'(B_{r})
\right] +E\left[f(B_{s})f''(B_{r})\right]\right\}\\
&\qquad+\left\{A_1(s,r,\varepsilon_i,j) -A_{21}(s,r,\varepsilon_i,j)
\right\} E\left[f''(B_{s})f(B_{r})\right]\\
&\qquad+\left\{A_3(s,r,\varepsilon_i,j) -A_{22}(s,r,\varepsilon_i,j)
\right\} E\left[f(B_{s})f''(B_{r})\right]
\end{align*}
with $i,j=1,2$ and $i\neq j$. By symmetry, we only need to show that
this holds for $i=1,j=2$. We will establish the
convergence~\eqref{Part-I-eq1} with $i=1,j=2$ in two steps.

{\bf Step A-1}. The following convergence hold:
\begin{align}\label{step2-eq1}
\frac1{\varepsilon_1^2\varepsilon_2}\int_0^tds\int_0^s A_{21}(s,r,
\varepsilon_1,2)\left\{ E\left[f'(B_{s})f'(B_{r})\right]
+E\left[f''(B_{s})f(B_{r})\right]\right\}dr \longrightarrow 0,\\
\label{step2-eq2}
\frac1{\varepsilon_1^2\varepsilon_2}\int_0^tds\int_0^s A_{22}(s,r,
\varepsilon_1,2)\left\{ E\left[f'(B_{s})f'(B_{r})\right]
+E\left[f(B_{s})f''(B_{r})\right]\right\}dr \longrightarrow 0,
\end{align}
as $\varepsilon_1,\varepsilon_2\to 0$. In order to prove the
convergence~\eqref{step2-eq1} we need to estimate
$$
A_{21}(s,r,\varepsilon_1,2).
$$
Notice that, by Lemma~\ref{lem2.4}
\begin{equation}\label{step1-eq1-2}
\begin{split}
\frac1{\varepsilon_1^2\varepsilon_2}|A_{21}(s,r,\varepsilon_1,2)|&\leq \frac1{\varepsilon_1^2\varepsilon_2}|E[B_{r}(B_{s+\varepsilon_1}-B_s)]|\\
&\qquad\cdot
\left(\left|\varepsilon_2E[B_{s}(B_{r+\varepsilon_1}-B_r)]\right|
+\left|\varepsilon_1E[B_{s}(B_{r+\varepsilon_2}-B_r)] \right|\right)\\
&\leq C_{H,K}\frac{r}{s}
\end{split}
\end{equation}
for $s>r>0$. We get
\begin{align*}
&\frac1{\varepsilon_1^2\varepsilon_2}\int_{\varepsilon_1}^t
\int_0^{s-\varepsilon_1} |A_{21}(s,r,
\varepsilon_1,2)|\left|E\left[f'(B_{s})f'(B_{r})\right]
+E\left[f''(B_{s})f(B_{r})\right]\right|drds\\
&\leq C_{H,K}\int_{\varepsilon_1}^t
\int_0^{s-\varepsilon_1}\frac{r}{s}
\left(\frac{1}{s^{3/4}r^{1/4}}+\frac{1}{
r^{3/4}s^{1/4}\sqrt{s-r}}\right)drds \int_{\mathbb{R}}
|f(x)|^2(s-r+x^2)\varphi_{s}(x)dx\\
&\leq C_{H,K}\|f\|_{\mathscr H}
\end{align*}
by Lemma~\ref{lem2.6}. Moreover, for
$\varepsilon_1<s<t,\;0<r<s-\varepsilon_1$ we have
\begin{equation}\label{Part-A-step1-eq1}
\begin{split}
\varepsilon_2E[B_{s}(B_{r+\varepsilon_1}-B_r)]&-
\varepsilon_1E[B_{s}(B_{r+\varepsilon_2}-B_r)]\\
&=2^{-K}\varepsilon_2 \left\{g_{s}(r+\varepsilon_1)-g_s(r)
-(s-r-\varepsilon_1) +(s-r)\right\}\\
&\qquad-2^{-K}\varepsilon_1 \left\{g_{s}(r+\varepsilon_2)
-g_s(r)-(s-r-\varepsilon_2) +(s-r)\right\}\\
&=2^{-K} \left\{[g_{s}(r+\varepsilon_1)-g_s(r)]\varepsilon_2
-[g_{s}(r+\varepsilon_2)-g_s(r)]\varepsilon_1\right\}\\
&=2^{-K} \left\{g'_{s}(\xi)-g'_{s}(\eta)\right\}\varepsilon_1
\varepsilon_2\\
&=2^{-K}
\left\{\frac{\xi^{2H-1}}{(s^{2H}+\xi^{2H})^{1-K}}
-\frac{\eta^{2H-1}}{(s^{2H}+\eta^{2H})^{1-K}}\right\}
\varepsilon_1\varepsilon_2
\end{split}
\end{equation}
for some $\xi\in (r,r+\varepsilon_1)$ and $\eta\in
(r,r+\varepsilon_2)$ by Mean Value Theorem, which implies that
\begin{align}\notag
\frac1{\varepsilon_1^2\varepsilon_2} |A_{21}(s,r,\varepsilon_1,2)|&=
|E[B_{r}(B_{s+\varepsilon_1}-B_s)]|\\  \notag
&\qquad\cdot
\left|\varepsilon_2E[B_{s}(B_{r+\varepsilon_1}-B_r)]
-\varepsilon_1E[B_{s}(B_{r+\varepsilon_2}-B_r)] \right|\\
\label{Part-A-step1-eq2} &\leq
C_{H,K}\frac{r}{s}\left|\frac{\xi^{2H-1}}{(s^{2H} +\xi^{2H})^{1-K}}
-\frac{\eta^{2H-1}}{(s^{2H}+\eta^{2H})^{1-K}} \right|\\
\notag
&\longrightarrow 0
\end{align}
for all $s>r>0$, as $\varepsilon_1,\varepsilon_2\to 0$. This proves
\begin{align*}
&\frac1{\varepsilon_1^2\varepsilon_2}\int_{\varepsilon_1}^t
\int_0^{s-\varepsilon_1}|A_{21}(s,r,
\varepsilon_1,2)|\left|E\left[f'(B_{s})f'(B_{r})\right]
+E\left[f''(B_{s})f(B_{r})\right]\right|dsdr\longrightarrow 0
\end{align*}
by Lebesgue's dominated convergence theorem. On the other hand,
Lemma~\ref{lem2.6} and~\eqref{step1-eq1-2} imply that
\begin{align*}
&\frac1{\varepsilon_1^2\varepsilon_2}\int_0^{\varepsilon_1}ds\int_0^s\frac{r}{s}|A_{21}(s,r,
\varepsilon_1,2)|\left|E\left[f'(B_{s})f'(B_{r})\right]
+E\left[f''(B_{s})f(B_{r})\right]\right|dr\\
&\leq C_{H,K}\int_0^{\varepsilon_1}ds\int_0^s\frac{r}{s}
\left(\frac{1}{s^{3/4}r^{1/4}}+\frac{1}{
r^{3/4}s^{1/4}\sqrt{s-r}}\right)dr\int_{\mathbb{R}}
|f(x)|^2(s+x^2)\varphi_{s}(x)dx\\
&=C_{H,K}\left(\int_0^{\varepsilon_1}(s+\sqrt{s})ds\int_{\mathbb{R}}
|f(x)|^2\varphi_{s}(x)dx+\int_0^{\varepsilon_1}
(1+\frac1{\sqrt{s}})ds\int_{\mathbb{R}}|f(x)x|^2
\varphi_{s}(x)dx\right)\\
&\leq C_{H,K}(\varepsilon_1+\sqrt{\varepsilon_1})\|f\|_{\mathscr
H}^2\longrightarrow 0\qquad (\varepsilon_1,\varepsilon_2\to 0)
\end{align*}
and
\begin{align*}
&\frac1{\varepsilon_1^2\varepsilon_2}
\int_{\varepsilon_1}^tds\int_{s-\varepsilon_1}^s |A_{21}(s,r,
\varepsilon_1,2)|\left|E\left[f'(B_{s})f'(B_{r})\right]
+E\left[f''(B_{s})f(B_{r})\right]\right|dsdr\\
&\leq C_{H,K}\int_{\varepsilon_1}^tds\int_{s-\varepsilon_1}^s
\frac{r}{s} \left(\frac{1}{s^{3/4}r^{1/4}}+\frac{1}{
r^{3/4}s^{1/4}\sqrt{s-r}}\right)dr\int_{\mathbb{R}}
|f(x)|^2(s+x^2)\varphi_{s}(x)dx\\
&= C_{H,K}\int_{\varepsilon_1}^t\left(s^{-3/4}(s^{\frac34+1}
-(s-\varepsilon_1)^{\frac34+1})+\sqrt{\varepsilon_1}\right)
ds\int_{\mathbb{R}} |f(x)|^2\varphi_{s}(x)dx\\
&\qquad+C_{H,K}\int_{\varepsilon_1}^t\frac1s
\left(s^{-3/4}(s^{\frac34+1}
-(s-\varepsilon_1)^{\frac34+1})+\sqrt{\varepsilon_1}\right)
ds\int_{\mathbb{R}} |f(x)x|^2\varphi_{s}(x)dx\\
&\leq C_{H,K}(\varepsilon_1+\sqrt{\varepsilon_1})\|f\|_{\mathscr
H}^2\longrightarrow 0\qquad (\varepsilon_1,\varepsilon_2\to 0).
\end{align*}
Thus, we have show that
\begin{align*}
\frac1{\varepsilon_1^2\varepsilon_2}&\int_0^t\int_0^s|A_{21}(s,r,
\varepsilon_1,2)|\left|E\left[f'(B_{s})f'(B_{r})\right]
+E\left[f''(B_{s})f(B_{r})\right]\right|drds\longrightarrow 0
\end{align*}
as $\varepsilon_1,\varepsilon_2\to 0$, which obtains the
convergence~\eqref{step2-eq1}. In a same way one can prove the
convergence~\eqref{step2-eq2}.

{\bf Step A-2}. The following convergence holds:
\begin{align}\label{step3-eq1}
\frac1{\varepsilon_1^2\varepsilon_2}\int_0^tds\int_0^s
\left\{A_1(s,r,\varepsilon_1,2) -A_{21}(s,r,\varepsilon_1,2)
\right\} E\left[f''(B_{s})f(B_{r})\right]dr \longrightarrow 0,\\
\label{step3-eq2}
\frac1{\varepsilon_1^2\varepsilon_2}\int_0^tds\int_0^s
\left\{A_3(s,r,\varepsilon_1,2) -A_{22}(s,r,\varepsilon_1,2)\right\}
E\left[f(B_{s})f''(B_{r})\right]dr \longrightarrow 0,
\end{align}
as $\varepsilon_1,\varepsilon_2\to 0$. We have
\begin{align}\notag
A_1(s,r,&\varepsilon_1,2)- A_{21}(s,r,\varepsilon_1,2)\\  \notag
&=\varepsilon_2 E[B_{s}(B_{r+\varepsilon_1}-B_r)]
E[B_{s}(B_{s+\varepsilon_1}-B_s)]\\  \notag
&\qquad-\varepsilon_1 E[B_{s}(B_{r+\varepsilon_2}-B_r)]
E[B_{s}(B_{s+\varepsilon_1}-B_s)]\\   \notag
&\qquad-\varepsilon_2
E[B_{s}(B_{r+\varepsilon_1}-B_r)]
E[B_{r}(B_{s+\varepsilon_1}-B_s)]\\    \notag
&\qquad+\varepsilon_1
E[B_{s}(B_{r+\varepsilon_2}-B_r)]
E[B_{r}(B_{s+\varepsilon_1}-B_s)]\\     \label{part3-eq010}
&=\left\{\varepsilon_2 E[B_{s}(B_{r+\varepsilon_1}-B_r)]-\varepsilon_1 E[B_{s}(B_{r+\varepsilon_2}-B_r)]\right\} E[(B_{s}-B_{r})(B_{s+\varepsilon_1}-B_s)],
\end{align}
which deduces
\begin{align*}
|A_1(s,r,\varepsilon_1,2)&- A_{21}(s,r,\varepsilon_1,2)|
\leq \varepsilon_2 |E[B_{s}(B_{r+\varepsilon_1}-B_r)]
E[(B_{s}-B_{r})(B_{s+\varepsilon_1}-B_s)]|\\
&\qquad+\varepsilon_1 |E[B_{s}(B_{r+\varepsilon_2}-B_r)]
E[(B_{s}-B_{r})(B_{s+\varepsilon_1}-B_s)]|\\
&\leq C_{H,K}\varepsilon_1^2\varepsilon_2\frac{s-r}{s}
\end{align*}
by Lemma~\ref{lem2.3}, Lemma~\ref{lem2.4} and the
estimate~\eqref{sec2-1-eq2.9=00}. It follows from Lemma~\ref{lem2.6}
that
\begin{align*}
\frac1{\varepsilon_1^2\varepsilon_2}\int_0^tds\int_0^s
|A_1(s,r,\varepsilon_1,2)- &A_{21}(s,r,\varepsilon_1,2)|
|E[f''(B_{s})f(B_{r})]|dr\\
&\leq \int_0^tds\int_0^s\frac{s-r}{s}|E[f''(B_{s})f(B_{r})]|dr\\
&\leq C_{H,K}\int_0^tds\int_0^s
\frac{1}{r^{1/4}s^{3/4}}dr\int_{\mathbb{R}}
f^2(x)\varphi_{s}(x)dx\\
&\leq C_{H,K}\|f\|_{\mathscr H}^2.
\end{align*}
On the other hand, by~\eqref{part3-eq010},~\eqref{Part-A-step1-eq1},
~\eqref{Part-A-step1-eq2} and Lemma~\ref{lem2.4} we have
\begin{align*}
\frac1{\varepsilon_1^2\varepsilon_2}& \left|A_1(s,r,
\varepsilon_1,2)- A_{21}(s,r,\varepsilon_1,2)\right|\longrightarrow
0,
\end{align*}
as $\varepsilon_1,\varepsilon_2\to 0$, for all $s>r>0$, which
implies that the convergence~\eqref{step3-eq1} holds by Lebesgue's
dominated convergence theorem. Similarly, one can
prove~\eqref{step3-eq2}.

{\bf Part B}. The following convergence hold:
\begin{align}\label{Part-B-eq1}
\frac1{\varepsilon_i^2\varepsilon_j}\int_0^t\int_0^r
\Upsilon(s,r,\varepsilon_i,j)dsdr \longrightarrow 0 \qquad
(i,j=1,2,i\neq j)
\end{align}
as $\varepsilon_1,\varepsilon_2\to 0$. For $r>s>0$ we can decompose
$\Upsilon(s,r,\varepsilon_i,j)$ as follows
\begin{align*}
\Upsilon(s,r,\varepsilon_i,j)
&=A_1(s,r,\varepsilon_i,j)E\left[f''(B_{s})f(B_{r})\right]
+A_3(s,r,\varepsilon_i,j)E\left[f(B_{s})f''(B_{r})\right]\\
&\qquad+\left(A_{21}(s,r,\varepsilon_i,j)+A_{22}(s,r,\varepsilon_i,j)
\right)E\left[f'(B_{s})f'(B_{r})\right]\\
&=A_{21}(s,r,\varepsilon_i,j)\left\{ E\left[f'(B_{s})f'(B_{r})
\right] +E\left[f(B_{s})f''(B_{r})\right]\right\}\\
&\qquad+A_{22}(s,r,\varepsilon_i,j)\left\{ E\left[f'(B_{s})f'(B_{r})
\right] +E\left[f''(B_{s})f(B_{r})\right]\right\}\\
&\qquad+\left\{A_3(s,r,\varepsilon_i,j) -A_{21}(s,r,\varepsilon_i,j)
\right\} E\left[f''(B_{s})f(B_{r})\right]\\
&\qquad+\left\{A_1(s,r,\varepsilon_i,j) -A_{22}(s,r,\varepsilon_i,j)
\right\} E\left[f(B_{s})f''(B_{r})\right]
\end{align*}
with $i,j=1,2$ and $i\neq j$. By the same method
proving~\eqref{step2-eq1} we can show that the following convergence
hold
\begin{align*}
\frac1{\varepsilon_i^2\varepsilon_j}\int_0^tdr\int_0^r
A_{21}(s,r,\varepsilon_i,j)\left\{ E\left[f'(B_{s})f'(B_{r})
\right] +E\left[f(B_{s})f''(B_{r})\right]\right\}ds
\longrightarrow 0,\\
\frac1{\varepsilon_i^2\varepsilon_j}\int_0^tdr\int_0^r A_{22}(s,r,
\varepsilon_i,j)\left\{ E\left[f'(B_{s})f'(B_{r})\right]
+E\left[f''(B_{s})f(B_{r})\right]\right\}ds \longrightarrow 0,
\end{align*}
as $\varepsilon_1,\varepsilon_2\to 0$. On the other hand, clearly,
we have
\begin{align*}
A_3(s,r,\varepsilon,j)&-A_{21}(s,r,\varepsilon,j)\\
&=\varepsilon_jE\left[B_{r}(B_{s+\varepsilon}-B_s)\right]
E\left[(B_{r}-B_{s})(B_{r+\varepsilon}-B_r)\right]\\
&\qquad-\varepsilon E\left[B_{r}(B_{s+\varepsilon_1}-B_s)\right]
E\left[(B_{r}-B_{s})(B_{r+\varepsilon_2}-B_r)\right]\\
A_1(s,r,\varepsilon,j)&-A_{22}(s,r,\varepsilon,j)\\
&=-\varepsilon_jE\left[B_{s}(B_{s+\varepsilon}-B_s)\right]
E\left[(B_{r}-B_{s})(B_{r+\varepsilon}-B_r)\right]\\
&\qquad+\varepsilon E\left[B_{s}(B_{s+\varepsilon_1}-B_s)\right]
E\left[(B_{r}-B_{s})(B_{r+\varepsilon_2}-B_r)\right]
\end{align*}
for all $r>s>0$. Thus, in the same way as proof of
~\eqref{step3-eq1} and~\eqref{step3-eq2} one can prove the
convergence
\begin{align*}
\frac1{\varepsilon_i^2\varepsilon_j}&\int_0^tdr\int_0^r
\left\{A_3(s,r,\varepsilon_i,j)-A_{21}(s,r,\varepsilon_i,j)
\right\} E\left[f''(B_{s})f(B_{r})\right]ds\longrightarrow 0\\
\frac1{\varepsilon_i^2\varepsilon_j}&\int_0^tdr\int_0^r
\left\{A_1(s,r,\varepsilon_i,j) -A_{22}(s,r,\varepsilon_i,j)
\right\} E\left[f(B_{s})f''(B_{r})\right]ds\longrightarrow 0
\end{align*}
as $\varepsilon_1,\varepsilon_2\to 0$, and the
convergence~\eqref{Part-B-eq1} follows.

{\bf Part C}. The following convergence holds:
\begin{align}\label{Part-C-eq1}
\frac1{\varepsilon_i^2\varepsilon_j}\int_0^tds\int_0^tA_4(s,r,
\varepsilon_i,j)E\left[f(B_{s})f(B_{r})\right]dr \longrightarrow
0\qquad (i,j=1,2,i\neq j)
\end{align}
as $\varepsilon_i,\varepsilon_j\to 0$. We have
\begin{align*}
A_4(s,r,\varepsilon_1,2)&=\varepsilon_2
E[(B_{r+\varepsilon_1}-B_r)(B_{s+\varepsilon_1}-B_s)]-\varepsilon_1
E[(B_{s+\varepsilon_1}-B_s)(B_{r+\varepsilon_2}-B_r)]\\
&=2^{-K}\Bigl(
[g_{s+\varepsilon_1}(r+\varepsilon_1)-g_{s}(r+\varepsilon_1)
-[g_{s+\varepsilon_1}(r)-g_{s}(r)]\varepsilon_2\\
&\qquad\qquad-[g_{s+\varepsilon_1}(r+\varepsilon_2)
-g_{s}(r+\varepsilon_2) -[g_{s+\varepsilon_1}(r)-g_{s}(r)]
\varepsilon_1\Bigr)\\
&\quad+2^{-K}\Bigl(\varepsilon_2\left[-|s-r|+|s+\varepsilon_1-r|
+|s-r-\varepsilon_1|-|s-r|\right]\\
&\qquad\qquad-\varepsilon_1\left[-|s+\varepsilon_1-r-\varepsilon_2|
+|s+\varepsilon_1-r| +|s-r-\varepsilon_2|-|s-r|\right]\Bigr)\\
&\equiv
2^{-K}A_{41}(s,r,\varepsilon_1,2)+2^{-K}A_{42}(s,r,\varepsilon_1,2)
\end{align*}
for $s,r>0$. By Mean Value Theorem we have
\begin{align}\label{Part-C-eq2}
A_{41}(s,r,\varepsilon_1,2)&
=\varepsilon_1\varepsilon_2\left([g'_{s+\varepsilon_1}(\xi)
-g'_{s}(\xi)]-[g'_{s+\varepsilon_1}(\eta) -g'_{s}(\eta)]\right)
\end{align}
for some $\xi\in (r,r+\varepsilon_1)$ and $\eta\in
(r,r+\varepsilon_2)$. Now, the convergence~\eqref{Part-C-eq1} will
be varied in three cases.

For $0<r,s<\varepsilon_1$. It is easy to verify that
\begin{align}\label{Part-C-eq3}
|A_{42}(s,r,\varepsilon_1,2)|&\leq 2\varepsilon_1\varepsilon_2.
\end{align}
Combining this with
\begin{align}\label{Part-C-eq4}
|g'_{y}(x)|=\frac{x^{2H-1}}{(y^{2H}+x^{2H})^{1-K}}
=\left(\frac{x^{2H}}{y^{2H}+x^{2H}}\right)^{1-K}\leq 1,\qquad
x,y\geq 0,
\end{align}
we get
\begin{align*}
\frac1{\varepsilon_1^2\varepsilon_2}\int_0^{\varepsilon_1}ds
\int_0^{\varepsilon_1}&|A_4(s,r,
\varepsilon_1,2)E\left[f(B_{s})f(B_{r})\right]|dr\\
&\leq \frac{2}{\varepsilon_1}\int_0^{\varepsilon_1}
\int_0^{\varepsilon_1}\left\{E\left[f^2(B_{s})\right]
+E\left[f^2(B_{r})\right]\right\}drds\\
&=\int_0^{\varepsilon_1}E\left[f^2(B_{s})\right]ds\longrightarrow 0,
\end{align*}
as $\varepsilon_1, \varepsilon_2\to 0$, by Lebesgue's dominated
convergence theorem. Similarly, we can show that the following
convergence holds:
\begin{align*}
\frac1{\varepsilon_1^2\varepsilon_2}\int_{\varepsilon_1}^tds
\int_{s-\varepsilon_1}^s&|A_4(s,r,
\varepsilon_1,2)E\left[f(B_{s})f(B_{r})\right]|dr\longrightarrow 0,
\end{align*}
as $\varepsilon_1, \varepsilon_2\to 0$. For $s>r+\varepsilon_1$, by
using Mean Value Theorem to the function
$$
x\mapsto g'_{s+\varepsilon_1}(x) -g'_{s}(x),\qquad x\geq 0,
$$
again, we get
\begin{align*}
A_{41}(s,r,\varepsilon_1,2)&
=\varepsilon_1\varepsilon_2(\xi-\eta)\left(
g''_{s+\varepsilon_1}(\theta) -g''_{s}(\theta)\right)
\end{align*}
for a $\theta\in (\xi\wedge\eta, \xi\vee\eta)$, which gives
\begin{align*}
|A_{41}(s,r,\varepsilon_1,2)|& \leq
\varepsilon_1^2\varepsilon_2\left|g''_{s+\varepsilon_1}(\theta)
-g''_{s}(\theta)\right|\longrightarrow 0
\end{align*}
for all $s>r>0$, as $\varepsilon_1,\varepsilon_2\to 0$. It follows
that
\begin{align*}
\frac1{\varepsilon_1^2\varepsilon_2}\int_{\varepsilon_1}^tds
\int_0^{s-\varepsilon_1}&|A_4(s,r,
\varepsilon_1,2)E\left[f(B_{s})f(B_{r})\right]|dr\longrightarrow 0,
\end{align*}
as $\varepsilon_1,\varepsilon_2\to 0$ because
$$
\frac1{\varepsilon_1^2\varepsilon_2}|A_4(s,r, \varepsilon_1,2)|\leq
\frac2s
$$
for $s>r+\varepsilon_1$. Finally, by symmetry we have that
\begin{align*}
\frac1{\varepsilon_1^2\varepsilon_2}\int_0^tdr\int_0^rA_4(s,r,
\varepsilon_1,2)E\left[f(B_{s})f(B_{r})\right]dr \longrightarrow 0,
\end{align*}
as $\varepsilon_1,\varepsilon_2\to 0$, and moreover, in the same way
we can establish the convergence
\begin{align*}
\frac1{\varepsilon_2^2\varepsilon_1}\int_0^tds\int_0^tA_4(s,r,
\varepsilon_2,1)E\left[f(B_{s})f(B_{r})\right]dr \longrightarrow 0,
\end{align*}
as $\varepsilon_1,\varepsilon_2\to 0$. Thus, we have established the
convergence~\eqref{sec4-Con-eq1}, and the theorem follows.
\end{proof}

\begin{corollary}
Let $2HK=1$. If $f$ is uniformly bounded, then the quadratic covariation $[f(B),B]$ exists in $L^2$ and
\begin{align}\label{sec3-eq3-2}
E\left|[f(B),B]_t\right|^2\leq \left(C_{H,K}\max_{x}|f(x)|\right)t^2
\end{align}
for all $0\leq t\leq T$.
\end{corollary}

\section{An It\^o formula}\label{sec4-1}
Our main object of this section is to explain and prove the
following theorem which gives a generalized It\^o formula.
\begin{theorem}\label{th4.2-}
Let $2HK=1$ and let $f\in {\mathscr H}$ be left continuous with right limits. If $F$ is an absolutely continuous function with the
derivative $F'=f$, then the following It\^o type formula holds:
\begin{equation}\label{sec5-eq5.3}
F(B)=F(0)+\int_0^tf(B_s)dB_s+2^{K-2}\left[f(B),B\right]_t.
\end{equation}
\end{theorem}
Clearly, the formula~\eqref{sec5-eq5.3} is an analogue of
F\"ollmer-Protter-Shiryayev's formula (see Eisenbaum~\cite{Eisen1},
F\"ollmer {\it et al}~\cite{Follmer}, Moret--Nualart~\cite{Moret},
Russo--Vallois~\cite{Russo2}, and the references therein). It is an
improvement in terms of the hypothesis on $f$ and it is also quite
interesting itself. As an application we get the relationship
between the forward (pathwise) integral and Skorohod integral
\begin{align*}
\int_0^tf(B_s)d^{-}B_s&=\int_0^tf(B_s)dB_s
+\frac12(2^{K-1}-1)\left[f(B),B\right]_t
\end{align*}
for all $f\in {\mathscr H}$ left continuous with right limits. The result weakens the hypothesis of differentiability for $f$ (see Russo-Tudor~\cite{Russo-Tudor})

Beside on the localization argument and smooth
approximation one can prove Theorem~\ref{th4.2-}. The so-called the localization argument is that one can localize the domain ${\rm Dom}(\delta^{H,K})$ of the operator $\delta^{H,K}$ (see Nualart~\cite{Nua2}). Suppose that $\{(\Omega_n, u_n), n\geq 1\}\subset {\mathscr F}\times {\rm Dom}(\delta^{H,K})$ is a localizing sequence for $u$, i.e., the sequence $\{(\Omega_n, u_n), n\geq 1\}$ satisfies
\begin{itemize}
\item [(i)] $\Omega_n\uparrow \Omega$, a.s.;
\item [(ii)] $u=u_n$ a.s. on $\Omega_n$.
\end{itemize}
If $\delta(u^{(n)})= \delta(u^{(m)})$ a.s. on $\Omega_n$ for all
$m\geq n$, then, the divergence $\delta^{H,K}$ is the random variable
determined by the conditions
$$
\delta^{H,K}(u)|_{\Omega_n}=\delta^{H,K} (u^{(n)})|_{\Omega_n}\qquad {\rm {
for\;\; all\;\;}} n\geq 1,
$$
but it may depend on the localizing sequence. Under the localization
argument {\em one may assume that the function $f\in {\mathscr H}$
is uniformly bounded}. In fact, for any $k\geq 0$ we may consider
the set
$$
\Omega_k=\left\{\sup_{0\leq t\leq T}|B_t|<k\right\}
$$
and let $f^{[k]}$ be a measurable function such that $f^{[k]}=f$ on
$[-k,k]$ and vanishes outside. Then $f^{[k]}$ is uniformly bounded
and $f^{[k]}\in {\mathscr H}$ for every $k\geq 0$. Set
$\frac{d}{dx}F^{[k]}=f^{[k]}$ and $F^{[k]}=F$ on $(-k,k)$. If the
formula~\eqref{sec5-eq5.3} is true for all uniformly bounded
functions, then we get the desired formula
$$
F^{[k]}(B_t)=F^{[k]}(0)+\int_0^t
f^{[k]}(B_s)dB_s+2^{K-2}\left[f^{[k]}(B),B\right]_t
$$
on the set $\Omega_k$. Letting $k$ tend to infinity we deduce the
It\^o formula~\eqref{sec5-eq5.3} for all $f\in {\mathscr H}$ being
left continuous with right limits. Thus, we may assume that {\em
$f\in {\mathscr H}$ is uniformly bounded in the next discussion}.
\begin{lemma}[Nualart~\cite{Nua2},
Es-sebaiy and Tudor~\cite{Es-sebaiy}]\label{complete2.1} Let
$\{u^{(n)}\}$ be a sequence such that $u_n\to u$ in $L^2$, as $n\to
\infty$ and let
$$
\delta^{H,K}(u^{(n)})=\int_0^Tu_s^{(n)}dB_s,\qquad n\geq 1
$$
exist in $L^2$. If $\delta^{H,K}(u^{(n)})\to G$ in $L^2$, then
$\delta^{H,K}(u)=\int_0^Tu_sdB_s$ exists in $L^2$ and equals to $G$.
\end{lemma}

\begin{lemma}\label{corollary3.1}
Let $f,f_1,f_2,\ldots\in {\mathscr H}$. If $f_n\to f$ in ${\mathscr
H}$ as $n$ tends to infinity, then we have
\begin{equation}\label{cor3.1-eq1}
\int_0^tf_n(B_s)d^{\pm}B_s\longrightarrow \int_0^tf(B_s)d^{\pm}B_s
\end{equation}
and
\begin{equation}\label{cor3.1-eq2}
[f_n(B),B]_t\longrightarrow [f(B),B]_t
\end{equation}
in $L^2$ as $n\to \infty$.
\end{lemma}
\begin{proof}
The lemma follows from
$$
E\left|\int_0^tf_n(B_s)d^{\pm}B_s-
\int_0^tf(B_s)d^{\pm}B_s\right|^2\leq C_{H,K}\|f_n-f\|_{\mathscr
H}^2\to 0,
$$
as $n$ tends to infinity.
\end{proof}
\begin{proof}[Proof of Theorem~\ref{th4.2-}]
If $F\in C^2({\mathbb R})$, this is It\^o's formula since
$$
[f(B),B]_t=\int_0^tf'(B_s)d[B,B]_s=2^{1-K}\int_0^tf'(B_s)ds.
$$
If $F\not\in C^2({\mathbb R})$, we let $F'=f\in {\mathscr H}$ be
uniformly bounded and left continuous. Consider the function $\zeta$
on ${\mathbb R}$ by
\begin{equation}
\zeta(x):=
\begin{cases}
ce^{\frac1{(x-1)^2-1}}, &{\text { $x\in (0,2)$}},\\
0, &{\text { otherwise}},
\end{cases}
\end{equation}
where $c$ is a normalizing constant such that $\int_{\mathbb
R}\zeta(x)dx=1$. Define the so-called mollifiers
\begin{equation}\label{sec4-eq00-4}
\zeta_n(x):=n\zeta(nx),\qquad n=1,2,\ldots
\end{equation}
and the sequence of smooth functions
\begin{align}\label{sec4-eq00-1}
F_n(x):=\int_{\mathbb
R}F(x-{y})\zeta_n(y)dy=\int_0^2F(x-\frac{y}n)\zeta(y)dy,\qquad
n=1,2,\ldots
\end{align}
for all $x\in \mathbb R$. Denote $f_n=F_n'$ for $n=1,2,\ldots$. Then
$F_n\in C^\infty({\mathbb R})$, $f_n\in C^{\infty}({\mathbb R})\cap
{\mathscr H}$ and
$$
f_n(x)=\int_{\mathbb R}f(x-{y})\zeta_n(y)dy
$$
for all $n\geq 1$. It is easy to check that $F_n',f_n,f_n'$ ($n\geq 1$) satisfy the condition~\eqref{sec2-Ito-con1} in
Theorem~\ref{theorem-Ito}. Hence, Skorohod integral $\int_0^tf_n(B_s)dB_s$ exists and It\^o's formula
\begin{equation}\label{sec3-eq3-Ito-1}
F_n(B_t)=F_n(0)+\int_0^tf_n(B_s)dB_s+\frac12\int_0^tf_n'(B_s)ds
\end{equation}
holds for all $n\geq 1$.

On the other hand, using Lebesgue's dominated
convergence theorem, one can prove that as $n$ tends to infinity, $f_n\to f$ in ${\mathscr H}$ and
$$
F_n(B_t)\longrightarrow F(B_t),\quad f_n(B_t)\longrightarrow f(B_t),
$$
in $L^2$, for all $t\in [0,T]$. Thus, we
get
\begin{align*}
2^{1-K}\int_0^tf_n'(B_s)ds=[f_n(B),B]_t \longrightarrow
\left[f(B),B\right]_t
\end{align*}
in $L^2$ by Lemma~\ref{corollary3.1}, as $n$ tends to infinity. It
follows that
\begin{align*}
\int_0^tf_n(B_s)dB_s&=F_n(B_t)-F_n(0)-
2^{K-2}[f_n(B),B]_t\\
&\longrightarrow F(B_t)-F(0)-2^{K-2}[f(B),B]_t
\end{align*}
in $L^2$, as $n$ tends to infinity. This completes the proof by
Lemma~\ref{complete2.1}.
\end{proof}


\section{The Bouleau-Yor identity}\label{sec4}
In this section we study one parameter integral of local time
$$
\int_{\mathbb R}f(x){\mathscr L}(dx,t),
$$
and establish the {\em Bouleau-Yor identity} between the integral
above and the quadratic covariation $[f(B),B]$, where $f$ is a
deterministic function and
$$
{\mathscr L}(x,t)=\int_0^t\delta(B_s-x)ds
$$
is the local time of bi-fBm $B$. Recall that the quadratic
covariation $[f(W),W]$ of Brownian motion $W$ can be characterized
as
$$
[f(W),W]_t=-\int_{\mathbb R}f(x){\mathscr L}^{W}(dx,t),
$$
where $f$ is locally square integrable and ${\mathscr L}^{W}(x,t)$ is the local time of Brownian motion $W$. This is called the {\em Bouleau-Yor identity}. More works for this can be found in
Bouleau-Yor~\cite{Bouleau}, Eisenbaum~\cite{Eisen1}, F\"ollmer {\it
et al}~\cite{Follmer}, Feng--Zhao~\cite{Feng},
Peskir~\cite{Peskir1}, Rogers--Walsh~\cite{Rogers2},
Yang--Yan~\cite{Yan2}, and the references therein. Moreover, this has be extended to fractional Brownian motion $B^H$ by Yan {\em et al}~\cite{Yan1,Yan2}.

\begin{lemma}\label{lem10.1}
Let $f_\triangle\in {\mathscr E}$. If
$$
f_\triangle=\sum_{j=1}^{N_1}x_j1_{(a_{j-1},a_j]} =\sum_{i=1}^{N_2}y_i1_{(b_{i-1},b_i]},
$$
we then have
\begin{align*}
\sum_jx_j[{\mathscr L}(a_j,t)-&{\mathscr
L}(a_{j-1},t)]=\sum_iy_i\left[{\mathscr L}(b_i,t)-{\mathscr L}(b_{i-1},t)\right]\\
&=-2^{K-1}[f_\triangle(B),B]_t.
\end{align*}
\end{lemma}
\begin{proof}
Take $F(x)=(x-a)^{+}-(x-b)^{+}$. Then $F$ is absolutely continuous with the derivative $F'=1_{(a,b]}\in {\mathscr H}$ being left continuous and bounded, and the It\^o formula~\eqref{sec5-eq5.3} yields
\begin{align*}
2^{K-2}\left[1_{(a,b]}(B),B\right]_t&=F(B_t)-F(0)-
\int_0^t1_{(a,b]}(B_s)dB_s
\end{align*}
for all $t\geq 0$. On the other hand, the Tanaka
formula~\eqref{eq2.4} follows
\begin{align*}
{\mathscr L}(a,t)-{\mathscr L}(b,t)=2F(B_t)-2F(0)-
2\int_0^t1_{(a,b]}(B_s)dB_s
\end{align*}
for all $t\geq 0$, which deduces
\begin{align*}
{\mathscr L}(a,t)-{\mathscr L}(b,t)=2^{K-1}\left[1_{(a,b]}(B),B\right]_t
\end{align*}
for all $t\geq 0$. Thus, the linearity property of the quadratic covariation implies that the lemma holds.
\end{proof}

As a direct consequence of Lemma~\ref{lem10.1} we can define the integral
\begin{equation}\label{sec5-eq5.1}
\int_{\mathbb
R}f_\triangle(x)\mathscr{L}(dx,t):=\sum_jx_j\left[{\mathscr L}(a_j,t)-{\mathscr L}(a_{j-1},t)\right]
\end{equation}
for every $f_\triangle\in {\mathscr E}$. Together this and Lemma~\ref{corollary3.1} lead to
\begin{equation*}
\lim_{n\to
\infty}\int_{\mathbb R}f_{\triangle,n}(x){\mathscr L}(dx,t)=\lim_{n\to
\infty}\int_{\mathbb R}\tilde{f}_{\triangle,n}(x){\mathscr L}(dx,t)\quad {\text { in $L^2$ }},
\end{equation*}
if $f_{\triangle,n}\to f$ and $\tilde{f}_{\triangle,n}\to f$
in ${{\mathscr H}}$, as $n$ tends to infinity, where
$\{f_{\triangle,n}\},\{\tilde{f}_{\triangle,n}\}\subset {\mathscr E}$. Thus, thanks to the density of ${\mathscr E}$ in ${{\mathscr H}}$, we can define integral of $f\in {\mathscr H}$  with respect to $x\mapsto {\mathscr L}(x,t)$ in the following manner:
\begin{equation*}
\int_{\mathbb R}f(x){\mathscr L}(dx,t):=\lim_{n\to
\infty}\int_{\mathbb R}f_{\triangle,n}(x){\mathscr L}(dx,t)\qquad {\text { in $L^2$ }},
\end{equation*}
provided $f_{\triangle,n}\to f$
in ${{\mathscr H}}$, as $n$ tends to infinity, where
$\{f_{\triangle,n}\}\subset {\mathscr E}$.

\begin{corollary}
Let $2HK=1$ and let $f\in {\mathscr H}$. Then the integral $\int_{\mathbb R}f(x){\mathscr L}(dx,t)$ exists in $L^2$, and
the Bouleau-Yor identity
\begin{equation}\label{sec5-eq5.2}
\left[f(B),B\right]_t=-2^{1-K}\int_{\mathbb R}f(x){\mathscr L}(dx,t)
\end{equation}
holds for all $t\in [0,T]$.
\end{corollary}

\begin{corollary}
Let $2HK=1$ and let $f,f_1,f_2,\ldots\in {\mathscr H}$. If $f_n\to f$ in ${\mathscr H}$, as $n$ tends to infinity, we then have
\begin{align*}
\int_{\mathbb R}f_n(x){\mathscr L}(dx,t)\longrightarrow
\int_{\mathbb R}f(x){\mathscr L}(dx,t)
\end{align*}
in $L^2$, as $n$ tends to infinity.
\end{corollary}

According to Theorem~\ref{th4.2-}, we get an analogue of
Bouleau-Yor's formula.
\begin{corollary}\label{cor6.1}
Let $2HK=1$ and let $f\in {\mathscr H}$ be left continuous with right limits. If $F$
is an absolutely continuous function with the derivative $F'=f$,
then the following It\^o type formula holds:
\begin{equation}
F(B_t)=F(0)+\int_0^tf(B_s)dB_s-\frac12\int_{\mathbb R}f(x){\mathscr
L}(dx,t).
\end{equation}
\end{corollary}
Recall that if $F$ is the difference of two convex functions, then $F$ is an absolutely continuous function with derivative of bounded variation. Thus, the It\^o-Tanaka formula (see Es-sebaiy and Tudor~\cite{Es-sebaiy})
\begin{align*}
F(B_t)&=F(0)+\int_0^tF^{'}(B_s)dB_s+\frac12\int_{\mathbb
R}{\mathscr L}(x,t)F''(dx)\\
&\equiv F(0)+\int_0^tF^{'}(B_s)dB_s-\frac12\int_{\mathbb
R}F'(x){\mathscr L}(dx,t)
\end{align*}
holds.


\begin{acknowledgement}
{\rm

The authors would like to thank the anonymous earnest referees whose
remarks and suggestions greatly improved the presentation of our
paper.

}
\end{acknowledgement}


\begin{thebibliography}{99}

\bibitem{Nua1}
E. Al\'os, O. Mazet and D. Nualart, Stochastic calculus with respect
to {G}aussian processes, {\it Ann. Probab.} {\bf 29} (2001),
766-801.

\bibitem{Bardina}
X. Bardina and C. Rovira, On It\^o formula for elliptic diffusion
processes, {\it Bernoulli}, {\bf 13} (2007), 820-830.

\bibitem{Bouleau}
N. Bouleau and M. Yor, Sur la variation quadratique des temps locaux
de certaines semimartingales, {\it C. R. Acad. Sci. Paris S\'er. I
Math.} {\bf 292} (1981), 491-494.

\bibitem{Eisen1}
N. Eisenbaum, Integration with respect to local time, { \it
Potential Anal.} {\bf 13} (2000), 303-328.

\bibitem{Eisen2}
N. Eisenbaum, Local time-space stochastic calculus for L\'evy
processes, {\it Stochastic Process. Appl.} {\bf 116} (2006),
757-778.

\bibitem{Elworthy}
K. D. Elworthy, A. Truman and H. Zhao, Generalized It\^o formulae
and space-time Lebesgue Stieltjes integrals of local times, {\it
S\'eminaire de Probabilit\'es} {\bf XL} (2007), 117-136.

\bibitem{Es-sebaiy}
K. Es-sebaiy and C. A. Tudor, Multidimensional bi-fractional
Brownian motion: It\^o and Tanaka formulas, {\it Stoch. Dyn.}, {\bf
7} (2007), 366-388.

\bibitem{Feng}
C. Feng and H. Zhao, Two-parameters $p,q$-variation paths and
integrations of local times, {\it Potential Anal.} {\bf 25} (2006),
165-204.

\bibitem{Follmer}
H. F\"ollmer, Ph. Protter and A. N. Shiryayev, Quadratic covariation
and an extension of It\^o's formula, {\it Bernoulli}, {\bf 1}
(1995), 149-169.

\bibitem{Geman}
D. Geman, J. Horowitz, Occupation densities, {\it Ann. Probab.} {\bf
8} (1980), 1-67.

\bibitem{Hou}
C. Houdr\'e and J. Villa,  An example of infinite dimensional
quasi-helix, {\it Stochastic models (Mexico City, 2002)},
pp.195-201, Contemp. Math., {\bf 336} (2003), Amer. Math. Soc.,
Providence, RI.

\bibitem{Jiang2}
Y. Jiang and Y. Wang, Self-intersection local times and collision
local times of bi-fractional Brownian motions, {\it Science in
China} {\bf Series A}: Mathematics, {\bf 52} (2009), 1905-1919.

\bibitem{Kruk}
I. Kruk, F. Russo and C. A. Tudor, Wiener integrals, Malliavin
calculus and covariance measure structure, {\it J. Funct. Anal.}
{\bf 249} (2007), 92-142.

\bibitem{Lei-Nualart}
P. Lei and D. Nualart, A decomposition of the bi-fractional Brownian
motion and some applications, {\em Statist. Probab. Lett.} {\bf 79}
(2009), 619-624.

\bibitem{Moret}
S. Moret and D. Nualart, Quadratic covariation and It\^o's formula
for smooth nondegenerate martingales, {\it J. Theoret. Probab.},
{\bf 13} (2000), 193-224.

\bibitem{Nua2}
D. Nualart, {\it Malliavin Calculus and Related Topics,} 2nd edition
Springer, New York (2006).

\bibitem{Peskir1}
G. Peskir, A change-of-variable formula with local time on curves,
{\it J. Theoret. Probab.} {\bf 18} (2005), 499-535.

\bibitem{Rogers2}
C. G. Rogers and J. B. Walsh, Local time and stochastic area
integrals, {\it Ann. Probab.} {\bf 19} (1991), 457-482.

\bibitem{Russo-Tudor}
F. Russo and C. A. Tudor, On the bi-fractional Brownian motion, {\it
Stochastic Process. Appl.} {\bf 5} (2006), 830--856.

\bibitem{Russo2}
F. Russo and P. Vallois, It\^o formula for ${\mathcal
C}^1$-functions of semimartingales, {\it Probab. Theory Rel.
Fields}. {\bf 104} (1996), 27-41.

\bibitem{Russo-Vallois2}
F. Russo and P. Vallois, Stochastic calculus with respect to a
continuous finite quadratic variation process, {\it Stochastics and
Stochastics Reports}, {\bf 70} (2000), 1-40.

\bibitem{Russo-Vallois3}
F. Russo and P. Vallois, Elements of stochastic calculus via
regularization, {\it S\'eminaire de Probabilit\'es} {\bf XL} (2007),
147-185.

\bibitem{Shen-yan}
G. Shen and L. Yan, Smoothness for the collision local times of
bifractional Brownian motions. {\it SCIENCE CHINA Math.} {\bf 54}
(2011), 1859-1873.

\bibitem{Tudor-Xiao}
C. A. Tudor and Y. Xiao, Some path properties of bi-fractional
brownian motion, {\it Bernoulli}, {\bf 13} (2007), 1023-1052.

\bibitem{Yan1}
L. Yan, C. Chen and J. Liu, The generalized quadratic covariation
for fractional Brownian motion with Hurst index less than $1/2$,
submitted (2011).

\bibitem{Yan4}
L. Yan, J. Liu and C. Chen, On the collision local time of
bi-fractional Brownian motions, {\it Stoch. Dyn.} {\bf 9} (2009),
479-491.

\bibitem{Yan2}
L.Yan, J. Liu and X. Yang, Integration with respect to fractional
local time with Hurst index $1/2<H<1$, {\it Potential Anal.}, {\bf
30} (2009), 115-138.

\bibitem{Yan3}
L.Yan and X. Yang, Some remarks on local time-space calculus, {\it
Statist. Probab. Lett.} {\bf 77} (2007), 1600-1607.

\end{thebibliography}
\end{document}